\documentclass[journal]{IEEEtran}      

\pdfoutput=1
\IEEEoverridecommandlockouts   

\newtheorem{Remark}{Remark}
\newtheorem{Corollary}{Corollary}

\newenvironment{Proof}{\noindent{\bf Proof.\/}}{\hfill $\blacksquare$\par}
\newtheorem{Theorem}{Theorem}
\newtheorem{Lemma}{Lemma}

\newtheorem{Assumption}{Assumption}

\newtheorem{Proposition}{Proposition}

\usepackage{graphicx}
\usepackage{amsmath,bm}
\usepackage{amssymb}
\usepackage{algorithm}
\usepackage{algpseudocode}
\usepackage{cite}
\usepackage{hyperref}
\usepackage[caption=false]{subfig}
\usepackage{color}
\hypersetup{hidelinks}

\allowdisplaybreaks

\algdef{SE}[DOWHILE]{Do}{doWhile}{\algorithmicdo}[1]{\algorithmicwhile\ #1}%

\newcommand{\mathactivatecomma}{%
  \begingroup\lccode`~=`\,
  \lowercase{\endgroup\edef~}{\mathchar\the\mathcode`\,\penalty0 }}

\algnewcommand{\Initialize}[1]{%
  \State \textbf{Initialize: $i \in \mathcal{V}$}
  \Statex \hspace*{\algorithmicindent}\parbox[t]{.8\linewidth}{\raggedright #1}
}

\algnewcommand{\Iteration}[1]{%
  \State \textbf{Iteration $(k\geq 0)$: $i \in \mathcal{V}$}
  \Statex \hspace*{\algorithmicindent}\parbox[t]{.8\linewidth}{\raggedright #1}
}

\algnewcommand{\Output}[1]{%
  \State \textbf{Output: $i \in \mathcal{V}$}
  \Statex \hspace*{\algorithmicindent}\parbox[t]{.8\linewidth}{\raggedright #1}
}

\title{\LARGE \bf
Gradient-Free Distributed Optimization with Exact Convergence
}

\author{Yipeng Pang and Guoqiang Hu
\thanks{This work was supported by Singapore Ministry of Education Academic Research Fund Tier 1 RG180/17(2017-T1-002-158).}
\thanks{Y. Pang and G. Hu are with the School of Electrical and Electronic Engineering, Nanyang
Technological University, 639798, Singapore
        {\tt\small ypang005@e.ntu.edu.sg, gqhu@ntu.edu.sg}.}%
}

\begin{document}

\bstctlcite{IEEEexample:BSTcontrol}

\maketitle
\thispagestyle{empty}
\pagestyle{empty}

\begin{abstract}
In this paper, a gradient-free distributed algorithm is introduced to solve a set constrained optimization problem under a directed communication network. Specifically, at each time-step, the agents locally compute a so-called pseudo-gradient to guide the updates of the decision variables, which can be applied in the fields where the gradient information is unknown, not available or non-existent. A surplus-based method is adopted to remove the doubly stochastic requirement on the weighting matrix, which enables the implementation of the algorithm in graphs having no associated doubly stochastic weighting matrix. For the convergence results, the proposed algorithm is able to obtain the exact convergence to the optimal value with any positive, non-summable and non-increasing step-sizes. Furthermore, when the step-size is also square-summable, the proposed algorithm is guaranteed to achieve the exact convergence to an optimal solution. In addition to the standard convergence analysis, the convergence rate of the proposed algorithm is also investigated. Finally, the effectiveness of the proposed algorithm is verified through numerical simulations.
\end{abstract}


\section{Introduction}
In recent years, with the prevalence of multi-agent systems, there has been a growing interest in solving the optimization problem in a distributed scheme. The advantage of doing so is that agents access local information and communicate with the neighbors only, making it suitable for the applications with large data size, huge computation and complex network structure, such as parameter estimation and detection \cite{Nowak2003,Ram2010}, source localization in sensor networks \cite{Lesser2003,Rabbat}, utility maximization \cite{Palomar2007}, resource allocation \cite{Chiang2007,Shen2012}, and multi-robot coordination \cite{Dong2017,Feng2017,Sun2018,Wu2019}. Distributed optimization of a sum of cost functions have been extensively studied over decades, such as the work in \cite{Chang2014,Nedic2010,Masubuchi2016,Zhu2012,Shi2015,Xu2017,Qu2017,Yuan2018,Lin2018,Yang2019}. A common underlying assumption in all these methods is that the derivative term of the local cost functions and the constraints can be directly accessed. However, there are many applications in the fields of bio-chemistry, aircraft design, hydro-dynamics, earth sciences, \textit{etc}., where the relation between the variables and the objective functions are unknown, the gradient information is not available for usage, or the derivative is not possible to determine \cite{Kramer2011}, these methods are no longer applicable. Hence, researchers start to draw attention to the gradient-free optimization.

Gradient-free optimization schemes can be traced back to the age of developing optimization theory, such as the work in \cite{Matyas1965}. Recent studies on this topic have been reported in \cite{Shamir2013,Nesterov2017,Yuan2015,Li2015,Chen2017,Yuan2015a,Pang2017,Pang2020}. Shamir \textit{et al.} in \cite{Shamir2013} investigated the performance of stochastic gradient descent method for non-smooth optimization problems. An averaging scheme was proposed to attain the minimax-optimal rates. On the other hand, Nesterov \textit{et al.} in \cite{Nesterov2017} provided an explicit way of computing the stochastic gradient information known as gradient-free oracle and investigated the convergence property for both convex and non-convex problems. This idea was extended to minimize a sum of non-smooth but Lipschitz continuous functions in \cite{Yuan2015,Li2015,Chen2017}, where the Gaussian smoothing technique was introduced to obtain the gradient-free oracle to replace the derivative in the standard subgradient methods. 
The same technique was applied to the algorithms in \cite{Yuan2015a} and \cite{Pang2017,Pang2020}, where the doubly stochastic requirement on the weighting matrix was removed by adopting a push-sum method \cite{Nedic2015} and a surplus-based method \cite{Cai2012}, respectively.
It should be noted that these derivative-free methods are based on the Gaussian smoothing technique, where the introduced smoothing parameter imposes an additional penalty term along the iteration. Thus, only an inexact convergence to a neighborhood of the optimal value can be achieved. To achieve the exact convergence, Duchi \textit{et al.} in \cite{Duchi2015} introduced a two point gradient estimation techique, and proved the exact convergence of the function value to the optimal value by choosing appropriate smoothing parameter sequences. This technique was extended to the distributed scenario in \cite{Yuan2016a,Pang2018} where an exact convergence of the function value to the optimal value was obtained.


In this paper, we aim to investigate gradient-free distributed optimization algorithms with exact convergence. Motivated by our work in \cite{Pang2018}, a distributed projected pseudo-gradient descent method is proposed to achieve an exact convergence with possibly a larger class of the step-sizes. The convergence properties of the proposed algorithm are carefully studied with different settings of the step-size. The main contributions of this work are summarized as follows.
\begin{enumerate}
\item Most gradient-free optimization algorithms, \textit{e.g.}, \cite{Nesterov2017,Yuan2015,Li2015,Chen2017,Yuan2015a,Pang2017,Pang2020} are based on Gaussian smoothing techniques, and hence can only achieve approximate convergence results. In terms of the exact convergence results, the work in \cite{Yuan2016a} proved an exact convergence of the function value to the optimal value for a step-size of $\alpha_k=\frac1{\sqrt{k}}$ ($k$ is the iteration index), and our work in \cite{Pang2018} proved the same convergence result for a non-summable and square-summable step-size. In this work, we introduce an optimal averaging scheme locally to trace a weighted average of the decision variable along the iteration. This averaging scheme is straightforward in terms of the implementation, and is able to obtain the exact convergence of the function value to the optimal value with any positive, non-increasing and non-summable step-sizes, hence increasing the range of the step-size selection. 
\item The convergence of the agent's function value does not imply that its decision variable also converges. The square-summable step-size condition is a typical setting in subgradient descent algorithms, \textit{e.g.}, \cite{Chang2014,Nedic2010,Masubuchi2016,Zhu2012,Xi2016,Lobel2011,Nedic2015,Mai2016,Xi2017a,Xiao2016,Pang2018} to establish the exact convergence of the agent's decision variable to an optimal solution. In this work, we show that this result also holds in distributed gradient-free algorithms. The proposed distributed projected pseudo-gradient descent method is guaranteed to achieve the exact convergence of the agent's decision variable to an optimal solution when the step-size also satisfies square-summable condition, which recovers the same convergence results in the literature.
\item The convergence rate has been widely studied in gradient-based distributed optimization literature, but received little attention in gradient-free distributed optimization literature. The only relevant works are \cite{Yuan2015a,Pang2020} and \cite{Li2015}, where \cite{Yuan2015a,Pang2020} proved a rate of $O(\frac{\ln t}{\sqrt{t}})$ for a diminishing step-size, and \cite{Li2015} showed a rate of $O(\frac1{\sqrt{t}})$ for a constant step-size if the number of iterations $t$ is known in advance. However, these rates were obtained for the algorithms with approximate convergence. In this work, the convergence rate of the proposed algorithm is studied, and we obtain the same convergence rate results as in \cite{Yuan2015a,Pang2020,Li2015} for the two settings of the step-size, but with exact convergence results.
\end{enumerate}

The rest of the paper is organized as follows. The problem is defined in Section~\ref{sec:problem_formulation}. Section~\ref{sec:distr_opt} introduces the proposed algorithm. The detailed convergence analysis is conducted in Section~\ref{sec:conv_analysis}, where some auxiliary lemmas are introduced, followed by the main results of the paper. The numerical simulations are presented in Section~\ref{sec:simulation} to illustrate the performance of the algorithm. Section~\ref{sec:conclusion} concludes the paper.

\section{Problem Formulation}\label{sec:problem_formulation}

For a directed graph $\mathcal{G} = \{\mathcal{V},\mathcal{E}\}$, $\mathcal{V} = \{1, 2, \ldots, N\}$ is the set of agents, and $\mathcal{E} \subset \mathcal{V}\times\mathcal{V}$ is the set of ordered pairs, ($i,j$), $i,j\in\mathcal{V}$, where agent $i$ is able to send information to agent $j$. We denote the set of agent $i$'s in-neighbors by $\mathcal{N}^{\text{in}}_i = \{j \in \mathcal{V} | (j,i)\in \mathcal{E}\}$ and out-neighbors by $\mathcal{N}^{\text{out}}_i = \{j \in \mathcal{V} | (i,j)\in \mathcal{E}\}$. Specifically, we allow both $\mathcal{N}^{\text{in}}_i$ and $\mathcal{N}^{\text{out}}_i$ to contain agent $i$ itself, and $\mathcal{N}^{\text{in}}_i\neq\mathcal{N}^{\text{out}}_i$ in general. The objective of the multi-agent system is to cooperatively solve the following set constrained optimization problem:
\begin{align}\label{eq:cost_function}
  \min f(\mathbf{x}) = \sum_{i=1}^N f_i(\mathbf{x}),\quad \mathbf{x} \in \mathcal{X},
\end{align}
where $\mathcal{X}\subseteq \mathbb{R}^n$ is a convex and closed set, and $f_i$ is a local cost function of agent $i$ and $\mathbf{x} = [x_1, \ldots, x_n]^\top$ is a global decision vector. The explicit expression of the local cost function $f_i$ is unknown, but the measurements can be made by agent $i$ only. Denote the (non-empty) solution set to \eqref{eq:cost_function} by $\mathcal{X}^\star$, \textit{i.e.}, $\mathcal{X}^\star = \arg\min_{\mathbf{x}\in\mathcal{X}} f(\mathbf{x})$.


Throughout this paper, we suppose the following assumptions hold:
\begin{Assumption}\label{assumption_graph}
The directed graph is strongly connected.
\end{Assumption}
\begin{Assumption}\label{assumption_local_f}
Each local cost function $f_i$ is convex, but not necessarily differentiable. For $\forall \mathbf{x} \in \mathcal{X}$, the subgradient $\partial f_i(\mathbf{x})$ exists and is bounded, \textit{i.e.}, there exists a positive constant $\hat{D}$ such that $\|\partial f_i(\mathbf{x})\|\leq \hat{D}$, $\forall \mathbf{x} \in \mathcal{X}$.
\end{Assumption}

\section{Algorithm}\label{sec:distr_opt}

In this section, we will develop the distributed projected pseudo-gradient descent method for the optimization problem defined in \eqref{eq:cost_function} as follows.

At time-step $k$, each agent $j$ broadcasts its state information $\mathbf{x}^j_k$ with a weighted auxiliary variable $[A_c]_{ij}\mathbf{y}^j_k$ to all of the nodes $i$ in its out-neighborhood. Then, for each agent $i$, on receiving the information $\mathbf{x}^j_k$, and $[A_c]_{ij}\mathbf{y}^j_k$ from all of the nodes $j$ in its in-neighborhood, it updates its variables $\mathbf{x}_{k+1}^i$ and $\mathbf{y}_{k+1}^i$\footnote{The update process does not require each agent to know the state information from its out-neighbors. but we assume agent $i$ knows the number of its in-neighbors and out-neighbors to design the weights in $A_r$ and $A_c$.}. Finally, each agent $i$ adopts an optimal averaging scheme to trace the average of $\mathbf{x}_{\ell}^i$, $\ell = 0,1,\ldots,k+1$ weighted by the step-size sequence, defined by $\widehat{\mathbf{x}}_{k+1}^i$. The updating law is given as follows.
\begin{subequations}\label{eq:algorithm}
\begin{align}
\mathbf{x}_{k+1}^i &= \mathcal{P}_{\mathcal{X}}\bigg[\sum_{j=1}^N [A_r]_{ij}\mathbf{x}^j_k + \epsilon \mathbf{y}^i_k - \alpha_k \mathbf{g}^i(\mathbf{x}^i_k)\bigg],\label{eq:update_x}\\
\mathbf{y}_{k+1}^i &= \mathbf{x}^i_k - \sum_{j=1}^N [A_r]_{ij}\mathbf{x}^j_k +  \sum_{j=1}^N [A_c]_{ij} \mathbf{y}^j_k - \epsilon \mathbf{y}^i_k,\label{eq:update_y}\\
\widehat{\mathbf{x}}_{k+1}^i &= \widehat{\mathbf{x}}_{k}^i + \frac{\alpha_{k+1}}{\sum_{\ell=0}^{k+1} \alpha_\ell}(\mathbf{x}_{k+1}^i-\widehat{\mathbf{x}}_{k}^i),\label{eq:update_xhat}
\end{align}
\end{subequations}
where $A_r,A_c$ are the row stochastic and column stochastic weighting matrices, respectively, \textit{i.e.,} $A_r\mathbf{1}_n = \mathbf{1}_n$, and $\mathbf{1}_n^\top A_c = \mathbf{1}_n^\top$. $\alpha_k > 0$ is a non-increasing step-size. $\epsilon$ is a small positive number. 
The auxiliary variable $\mathbf{y}^i_k$ is used to offset the shift caused by the unbalanced (non-doubly stochastic) weighting matrices ($A_r, A_c$), known as ``surplus''. The parameter $\epsilon$ is to specify the amount of surplus during the update (see \cite{Cai2012} for the details).
$\mathbf{g}^i(\mathbf{x}^i_k)$ is a pseudo-gradient motivated from \cite{Duchi2015}, given as
\begin{align}
\mathbf{g}^i(\mathbf{x}^i_k) &= \frac1{\beta_{2,k}}[f_i(\mathbf{x}^i_k+\beta_{1,k}\xi^i_{1,k}+\beta_{2,k}\xi^i_{2,k})\nonumber\\
&\qquad\qquad\qquad-f_i(\mathbf{x}^i_k+\beta_{1,k}\xi^i_{1,k})]\xi^i_{2,k},\label{grad_oracle}
\end{align}
$\beta_{1,k}$, $\beta_{2,k}$ are two positive non-increasing sequences with their ratio defined as 
\begin{align}\label{eq:ratio_two_sequence}
\tilde{\beta}_k = \beta_{2,k}/\beta_{1,k}.
\end{align}
$\xi^i_{1,k}$ and $\xi^i_{2,k} \in \mathbb{R}^n$ are two random variables satisfying the following assumption:

\begin{Assumption}\label{assumption_random_variables}
(Assumption F in \cite{Duchi2015}) The random variables $\xi^i_{1,k}$ and $\xi^i_{2,k} \in \mathbb{R}^n$ are generated by any one of the following: (a) both $\xi^i_{1,k}$ and $\xi^i_{2,k}$ are standard normal in $\mathbb{R}^n$ with identity covariance; (b) both $\xi^i_{1,k}$ and $\xi^i_{2,k}$ are uniform on the $\ell_2$-ball of radius $\sqrt{n+2}$; (c) the distribution of $\xi^i_{1,k}$ is uniform on the $\ell_2$-ball of radius $\sqrt{n+2}$ and the distribution of $\xi^i_{2,k}$ is uniform on the $\ell_2$-ball of radius $\sqrt{n}$.
\end{Assumption}

Similar to the gradient-free oracle in \cite{Nesterov2017}, at each time $k$, the pseudo-gradient operator \eqref{grad_oracle} estimates the gradient in a random direction $\xi^i_{2,k}$ with a parameter $\beta_{2,k}$, but the function difference is taken at a perturbed point $\mathbf{x}^i_k+\beta_{1,k}\xi^i_{1,k}$ instead of $\mathbf{x}^i_k$, where the amount of perturbation is determined by the parameter $\beta_{1,k}$ and the random variable $\xi^i_{1,k}$. As compared to the gradient-free oracle where the function difference is evaluated at $\mathbf{x}^i_k$ which may not be differentiable for non-smooth problems, the extra perturbation step in pseudo-gradient operator allows the function difference to be evaluated at a point which is less likely to be non-smooth. 
In fact, we can define a smoothed function of $f_i(\mathbf{x})$ based on the convolution of this perturbation, given by \cite{Duchi2015}, 
\begin{align*}
  f_{i,\beta_{1,k}}(\mathbf{x}) &= \mathbb{E}[f(\mathbf{x}+\beta_{1,k}\xi^i_{1,k})]\\
  & = \int_{\mathbb{R}^n} f_i(\mathbf{x}+\beta_{1,k}\xi^i_{1,k})d\mu(\xi^i_{1,k}),
\end{align*}
with the random variable $\xi^i_{1,k}\in\mathbb{R}^n$ having density $\mu$ with respect to Lebesgue measure.
$\beta_{1,k}$ is a positive non-increasing sequence.




In fact, algorithms \eqref{eq:update_x} and \eqref{eq:update_y} can be written into an equivalent form
\begin{align}\label{eq:d-dgd}
\mathbf{z}^i_{k+1} = \sum_{j=1}^{2N}[A]_{ij}\mathbf{z}^j_k + g^i_k,
\end{align}
where $g^i_k$, $i\in\{1,\ldots,2N\}$ is an augmented pseudo-gradient defined by $g^i_k = \mathbf{x}^i_{k+1}-\sum_{j=1}^N [A_r]_{ij}\mathbf{x}^j_k - \epsilon \mathbf{y}^i_k$ for $i\in\{1,\ldots,N\}$, $g^i_k = \mathbf{0}_n$ for $i\in\{N+1,\ldots,2N\}$; matrix $A\in\mathbb{R}^{2N\times2N}$ is an augmented weighting matrix defined by $A = [\begin{smallmatrix} A_r & \epsilon I \\ I-A_r & A_c - \epsilon I\end{smallmatrix}]$; and decision variable $\mathbf{z}^i_k$, $i\in\{1,\ldots,2N\}$ is defined by $\mathbf{z}^i_k = \mathbf{x}^i_k$ for $i\in\{1,\ldots,N\}$, $\mathbf{z}^i_k = \mathbf{y}^{i-N}_k$ for $i\in\{N+1,\ldots,2N\}$.

\section{Convergence Analysis}\label{sec:conv_analysis}

In this section, the detailed convergence analysis of our proposed algorithm is provided. We first introduce some auxiliary lemmas in Subsection~\ref{subsec:auxiliary_lemma}, followed by the main results in Subsection~\ref{subsec:main_results}.


\subsection{Auxiliary Lemmas}\label{subsec:auxiliary_lemma}

In this part, we introduce some auxiliary results, which will be helpful in the analysis of the main theorems.
We denote the $\sigma$-field generated by the entire history of the random variables from step 0 to $k-1$ by $\mathcal{F}_k$, \textit{i.e.,} $\mathcal{F}_k = \{(\mathbf{x}^i_0, i \in \mathcal{V});(\xi^i_{1,s}, \xi^i_{2,s}, i \in \mathcal{V}); 0\leq s\leq k-1\}$ with $\mathcal{F}_0 = \{\mathbf{x}^i_0, i \in \mathcal{V}\}$.

The following lemma summarizes some properties of function $f_{i,\beta_{1,k}}(\mathbf{x})$ and the pseudo-gradient $\mathbf{g}^i(\mathbf{x}^i_k)$.

\begin{Lemma}\label{lemma:property_f_mu}
(see \cite{Duchi2015}) Suppose Assumptions~\ref{assumption_local_f} and \ref{assumption_random_variables} hold. Then, for each $i \in \mathcal{V}$, the following properties of the function $f_{i,\beta_{1,k}}(\mathbf{x})$ are satisfied:
\begin{enumerate}
\item $f_{i,\beta_{1,k}}(\mathbf{x})$ is convex and differentiable, and it satisfies
\begin{align*}
f_i(\mathbf{x})\leq f_{i,\beta_{1,k}}(\mathbf{x})\leq f_i(\mathbf{x}) + \beta_{1,k}\hat{D}\sqrt{n+2},
\end{align*}
\item the pseudo-gradient $\mathbf{g}^i(\mathbf{x}^i_k)$ satisfies
\begin{align*}
\mathbb{E}[\mathbf{g}^i(\mathbf{x}^i_k)|\mathcal{F}_k] = \nabla f_{i,\beta_{1,k}}(\mathbf{x}^i_k)+\tilde{\beta}_k\hat{D}\mathbf{v},
\end{align*}
\item there is a universal constant $Q$ such that
\begin{align*}
\mathbb{E}[\|\mathbf{g}^i(\mathbf{x}^i_k)\||\mathcal{F}_k] &\leq \sqrt{\mathbb{E}[\|\mathbf{g}^i(\mathbf{x}^i_k)\|^2|\mathcal{F}_k]}\leq Q\mathcal{T}_k,
\end{align*}
\end{enumerate}
where $\beta_{1,k}$ and $\tilde{\beta}_k$ are defined in \eqref{eq:ratio_two_sequence}, $\mathbf{v} \in \mathbb{R}^n$ is a vector satisfying $\|\mathbf{v}\| \leq n\sqrt{3n}/2$, and $\mathcal{T}_k = \hat{D}\sqrt{n\Big[n\sqrt{\tilde{\beta}_k}+1+\ln n\Big]}$. If $\tilde{\beta}_k$ is bounded, then $\mathcal{T}_k$ is bounded. In this case, we denote the upper bound of $Q\mathcal{T}_k$ by $\mathcal{K}_1$.
\end{Lemma}

Following the results in \cite{Cai2012,Xi2016}, we have the following lemma on the convergence of the augmented weighting matrix $A$ in \eqref{eq:d-dgd}.
\begin{Lemma}\label{lemma:A_matrix}
Suppose Assumption~\ref{assumption_graph} holds. Let $\epsilon$ be the constant in the augmented weighting matrix $A$ in \eqref{eq:d-dgd} such that $\epsilon \in (0,\bar{\epsilon})$ with $\bar{\epsilon} = (\frac{1-|\lambda_3|}{20+8N})^N$, where $\lambda_3$ is the third largest eigenvalue of matrix $A$ with $\epsilon =0$. Then $\forall i,j \in \{1,\ldots,2N\}$, the entries $[A^k]_{ij}$ converge to their limits as $k \to \infty$ at a geometric rate, \text{i.e.,}
\begin{equation*}
\left\|A^k - \begin{bmatrix} \frac{\mathbf{1}_N\mathbf{1}^T_N}N &\frac{\mathbf{1}_N\mathbf{1}^T_N}N \\ \mathbf{0}_{N\times N} & \mathbf{0}_{N\times N} \end{bmatrix}\right\|_\infty \leq \Gamma \gamma^k, \quad k\geq 1,
\end{equation*}
where $\Gamma > 0$ and 
\begin{align*}
\gamma = \max\{|\lambda_3| + (20+8N)\epsilon^{\frac1N}, |\lambda_2(\epsilon)|\} \in(0,1)
\end{align*}
are some constants, and $\lambda_2(\epsilon)$ is the eigenvalue of the weighting matrix $A$ corresponding to the second largest eigenvalue $\lambda_2$ of matrix $A$ with $\epsilon =0$. 
\end{Lemma}
\begin{Proof}
The first part of the result follows directly from the proof of Lemma~1 in \cite{Xi2016}, where constant $\gamma$ is determined by the magnitude of the second largest eigenvalue of matrix $A$. Next we aim to characterize the second largest eigenvalue of matrix $A$. To do so, we denote $A$ by $A(\epsilon)$ to represent the dependency of $A$ on parameter $\epsilon$. Then, matrix $A(\epsilon)$ can be viewed as matrix $A(0)$ with some perturbations on $\epsilon$, where matrix $A(0)$ is matrix $A(\epsilon)$ by setting $\epsilon = 0$. Denote the eigenvalues of matrix $A(0)$ by $\lambda_1$, $\lambda_2$,...,$\lambda_{2N}$ with $|\lambda_1|\geq\cdots\geq|\lambda_{2N}|$. From the proof of Theorem~4 in \cite{Cai2012}, it holds that $1 = \lambda_1 = \lambda_2 >|\lambda_3|\geq\cdots\geq|\lambda_{2N}|$. After perturbation, we denote by $\lambda_i(\epsilon)$ the eigenvalues of matrix $A(\epsilon)$ corresponding to $\lambda_i$, $i=\{1,\ldots,N\}$. It should be noted that the eigenvalues of the perturbed matrix $A(\epsilon)$ do not necesssarily satisfy $|\lambda_1(\epsilon)|\geq\cdots\geq|\lambda_{2N}(\epsilon)|$ given that $|\lambda_1|\geq\cdots\geq|\lambda_{2N}|$. From Lemmas~10 and 11 in \cite{Cai2012}, when $\epsilon\in(0,\bar{\epsilon})$, we have the following inequality characterizing the distance between the corresponding eigenvalues $\lambda_i(\epsilon)$ and $\lambda_i$, $i=\{1,\ldots,N\}$
\begin{align*}
|\lambda_i(\epsilon)-\lambda_i| < 4(4+2N+\epsilon)\epsilon^{\frac1N}<(20+8N)\epsilon^{\frac1N}, 
\end{align*}
which gives $|\lambda_i(\epsilon)| < |\lambda_i| + (20+8N)\epsilon^{\frac1N}$. Hence, for $i=\{3,\ldots,N\}$, the above inequality yields $|\lambda_i(\epsilon)| < |\lambda_3| + (20+8N)\epsilon^{\frac1N}<1$ since $|\lambda_3|\geq\cdots\geq|\lambda_{2N}|$ and $\epsilon\in(0,\bar{\epsilon})$. Moreover, from the proof of Theorem~4 and Lemma~12 in \cite{Cai2012}, when $\epsilon\in(0,\bar{\epsilon})$, we have $\lambda_1(\epsilon) = 1$ and $|\lambda_2(\epsilon)| < 1$. Hence $\gamma$ can be selected as $\max\{|\lambda_3| + (20+8N)\epsilon^{\frac1N}, |\lambda_2(\epsilon)|\}$, which completes the proof.
\end{Proof}

\begin{Remark}\label{remark:effect_on_gamma}
The work in \cite{Xiao2004} has proposed solutions on the design of the weighting matrix to guarantee the fastest averaging speed when the weighting matrix is symmetric and doubly-stochastic. For the weighting matrix $A$ in this work, Lemma~\ref{lemma:A_matrix} shows that the averaging speed depends on constant $\gamma$.
From the proof of Lemma~\ref{lemma:A_matrix}, we can infer the effects of parameter $\epsilon$, the communication topology, and the number of agents $N$ on constant $\gamma$. For the effect of parameter $\epsilon$, noting that $|\lambda_2(\epsilon)|= 1$ when $\epsilon = 0$ and $|\lambda_3| + (20+8N)\epsilon^{\frac1N} = 1$ when $\epsilon=\bar{\epsilon}$, hence $\gamma$ is dominant by $|\lambda_2(\epsilon)|$ when $\epsilon$ is small, and then dominant by $|\lambda_3| + (20+8N)\epsilon^{\frac1N}$ when $\epsilon$ is large. That implies there is an optimal value of $\epsilon$ such that $\gamma$ is minimized (when $|\lambda_2(\epsilon)|=|\lambda_3| + (20+8N)\epsilon^{\frac1N}$). For the effect of the communication topology, suppose $\epsilon$ is set at the optimal value, then a graph with a smaller $|\lambda_3|$ leads to a smaller $\gamma$. For the effect of the number of agents $N$, since $|\lambda_3| + (20+8N)\epsilon^{\frac1N}$ is smaller for a smaller $N$, hence $\gamma$ is smaller for a smaller number of agents.
\end{Remark}


Define $\bar{\mathbf{z}}_k = \frac1N\sum_{i=1}^{2N}\mathbf{z}^i_k = \frac1N\sum_{i=1}^N\mathbf{x}^i_k + \frac1N\sum_{i=1}^N\mathbf{y}^i_k$, which is an average of $\mathbf{x}^i_k+\mathbf{y}^i_k $ over all agents at time-step $k$; and 
\begin{equation}\label{eq:def_zhat}
\widehat{\mathbf{z}}_{k}=\frac{\sum_{\ell=0}^k \alpha_\ell\bar{\mathbf{z}}_\ell}{\sum_{\ell=0}^k \alpha_\ell}, 
\end{equation}
which is an average of $\bar{\mathbf{z}}$ weighted by the step-size sequence $\alpha_\ell$ over time duration $k$.
Then, we can quantify the bounds of the terms $\mathbf{x}^i_k - \bar{\mathbf{z}}_k$ and $\mathbf{y}^i_k-\mathbf{0}_n$ as shown in the following lemma. For easy representation, we denote the aggregated norm of the augmented pseudo-gradient $\sum_{j=1}^{N}\|g^j_k\|$ by $\bm{G}_k$ in the rest of the paper.

\begin{Lemma}\label{lemma:bound_on_consensus}
Suppose Assumptions~\ref{assumption_graph}, \ref{assumption_local_f} and \ref{assumption_random_variables} hold. Let $\epsilon$ be the constant such that $\epsilon \in (0,\bar{\epsilon})$, where $\bar{\epsilon}$ is defined in Lemma~\ref{lemma:A_matrix}. Let $\{\mathbf{x}^i_k\}_{k\geq0}$ and $\{\mathbf{y}^i_k\}_{k\geq0}$ be the sequences generated by \eqref{eq:update_x} and \eqref{eq:update_y}, respectively. Then, it holds that
for $k \geq 1$
\begin{enumerate}
\item
$
\!
\begin{aligned}[t]
\|\mathbf{x}^i_k-\bar{\mathbf{z}}_k\|&\leq 2N\varsigma\Gamma\gamma^k+\Gamma\sum_{r=1}^{k-1}\gamma^{k-r}\bm{G}_{r-1}+ \bm{G}_{k-1};
\end{aligned}
$
\item
$
\!
\begin{aligned}[t]
\|\mathbf{y}^i_k\| &\leq 2N\varsigma\Gamma\gamma^k+\Gamma\sum_{r=1}^{k-1}\gamma^{k-r}\bm{G}_{r-1},
\end{aligned}
$
\end{enumerate}
where $\varsigma = \max\{\|\mathbf{x}^i_0\|,\|\mathbf{y}^i_0\|,i\in\mathcal{V}\}$, $\Gamma$ and $\gamma$ are the constants defined in Lemma~\ref{lemma:A_matrix}.
\end{Lemma}

\begin{Proof}
For $k \geq 1$, we have 
\begin{align}\label{eq:z_1}
\mathbf{z}^i_k = \sum_{j=1}^{2N}[A^k]_{ij}\mathbf{z}^j_0 + \sum_{r=1}^{k-1}\sum_{j=1}^{2N}[A^{k-r}]_{ij}g^j_{r-1} + g^i_{k-1}.
\end{align}
by applying (\ref{eq:d-dgd}) recursively. Then we can obtain that
\begin{equation}\label{eq:z_bar_2}
\begin{aligned}
\bar{\mathbf{z}}_k = \frac1N\sum_{j=1}^{2N}\mathbf{z}^j_0 + \frac1N\sum_{r=1}^{k-1}\sum_{j=1}^{2N}g^j_{r-1} + \frac1N\sum_{j=1}^{2N}g^j_{k-1},
\end{aligned}
\end{equation}
where we used column stochastic property of $A$.

For part (1), subtracting (\ref{eq:z_bar_2}) from (\ref{eq:z_1}) and taking the norm, we have that for $1 \leq i \leq N$ and $k\geq 1$,
\begin{align}
&\|\mathbf{z}^i_k-\bar{\mathbf{z}}_k\|\leq \sum_{j=1}^{2N}\bigg|[A^k]_{ij}-\frac1N\bigg|\varsigma+ \sum_{r=1}^{k-1}\sum_{j=1}^{N}\bigg|[A^{k-r}]_{ij}\nonumber\\
&\quad-\frac1N\bigg|\|g^j_{r-1}\|+ \frac{N-1}N\|g^i_{k-1}\| + \frac1N\sum_{j\neq i}\|g^j_{k-1}\|.\label{eq:z_bar_minus_z_2}
\end{align}
Noting that
$\frac{N-1}N\|g^i_{k-1}\| + \frac1N\sum_{j\neq i}\|g^j_{k-1}\|\leq\bm{G}_{k-1}$,
and applying the property of $[A^k]_{ij}$ from Lemma~\ref{lemma:A_matrix} to (\ref{eq:z_bar_minus_z_2}), we complete the proof of part (1).

For part (2), taking the norm in (\ref{eq:z_1}) for $N+1 \leq i \leq 2N$ and $k>1$, and applying the property of $[A^k]_{ij}$ from Lemma~\ref{lemma:A_matrix}, we complete the proof of part (2).
\end{Proof}

It can be seen from Lemma~\ref{lemma:bound_on_consensus} that the bound for the consensus terms is a function of the aggregated norm of the augmented pseudo-gradient term $\bm{G}_k$. Hence, in the following lemma, we provide some properties on this term $\bm{G}_k$.

\begin{Lemma}\label{lemma:bound_on_g}
Suppose Assumptions~\ref{assumption_graph}, \ref{assumption_local_f} and \ref{assumption_random_variables} hold. Let $\epsilon$ be the constant such that $0<\epsilon < \min(\bar{\epsilon},\frac{1-\gamma}{2\sqrt{3}N\Gamma\gamma})$, where $\bar{\epsilon}$, $\Gamma$ and $\gamma$ are the constants defined in Lemma~\ref{lemma:A_matrix}. Let $\tilde{\beta}_k$ defined in \eqref{eq:ratio_two_sequence} be bounded. Then, for any $K\geq 1$, the aggregated norm of the augmented pseudo-gradient term $\bm{G}_k$ satisfies that
\begin{enumerate}
\item
$
\!
\begin{aligned}[t]
\sum_{k=1}^K \alpha_k\mathbb{E}[\bm{G}_k] \leq \Phi_1 \sum_{k=1}^K \alpha^2_k + \Psi_1,
\end{aligned}
$
\item
$
\!
\begin{aligned}[t]
\sum_{k=1}^K \mathbb{E}[\bm{G}^2_k] \leq \Phi_2 \sum_{k=1}^K \alpha^2_k + \Psi_2,
\end{aligned}
$
\item
$
\!
\begin{aligned}[t]
\sum_{k=1}^K \sum_{i=1}^N\alpha_k\mathbb{E}[\|\mathbf{g}^i(\mathbf{x}^i_k)\|\bm{G}_k] \leq \Phi_3 \sum_{k=1}^K \alpha^2_k + \Psi_3,
\end{aligned}
$
\end{enumerate}
where $\Phi_1$, $\Psi_1$, $\Phi_2$, $\Psi_2$, $\Phi_3$ and $\Psi_3$ are positive bounded constants, and $\alpha_k>0$ is a non-increasing step-size.
\end{Lemma}
\begin{Proof} See Appendix A. \end{Proof}

In addition, we will frequently utilize the Stolz-Cesaro Theorem\cite{Nagy} to facilitate the analysis, which is quoted below for completeness.
\begin{Lemma}\label{theorem:Stolz-Cesaro}
(Stolz-Cesaro Theorem) If $\{b_k\}_{k\geq1}$ is a sequence of positive real numbers, such that $\sum_{k=1}^\infty b_k = \infty$, then for any sequence $\{a_k\}_{k\geq1}$ one has the inequality:
\begin{align*}
\liminf_{k\to\infty}\frac{a_k}{b_k}&\leq\liminf_{k\to\infty}\frac{a_1+a_2+\cdots+a_k}{b_1+b_2+\cdots+b_k}\\
&\leq\limsup_{k\to\infty}\frac{a_1+a_2+\cdots+a_k}{b_1+b_2+\cdots+b_k}\leq \limsup_{k\to\infty}\frac{a_k}{b_k}.
\end{align*}
In particular, if the sequence $\{a_k/b_k\}_{k\geq1}$ has a limit, then
\begin{align*}
\lim_{k\to\infty}\frac{a_1+a_2+\cdots+a_k}{b_1+b_2+\cdots+b_k}= \lim_{k\to\infty}\frac{a_k}{b_k}.
\end{align*}
\end{Lemma}

With the above lemmas, we are able to establish a one-step iteration and a consensus result under only non-summable step-size condition.

\begin{Proposition}\label{prop:consensus}
Suppose Assumptions~\ref{assumption_graph}, \ref{assumption_local_f} and \ref{assumption_random_variables} hold. Let $\{\mathbf{x}^i_k\}_{k\geq0}$, $\{\mathbf{y}^i_k\}_{k\geq0}$ and $\{\widehat{\mathbf{x}}^i_k\}_{k\geq0}$ be the sequences generated by \eqref{eq:algorithm} with a non-increasing step-size sequence $\{\alpha_k\}_{k\geq0}$ satisfying 
\begin{align*}
\sum_{k=0}^\infty\alpha_k= \infty, \quad\lim_{k\to\infty} \alpha_k=\alpha_\infty.
\end{align*} 
Let $\epsilon$ be the constant such that $0<\epsilon < \min(\bar{\epsilon},\frac{1-\gamma}{2\sqrt{3}N\Gamma\gamma})$, where $\bar{\epsilon}$, $\Gamma$ and $\gamma$ are the constants defined in Lemma~\ref{lemma:A_matrix}. Let $\tilde{\beta}_k$ defined in \eqref{eq:ratio_two_sequence} be bounded. Then

(1) $\widehat{\mathbf{x}}^i_k$ holds that
\begin{align*}
&\mathbb{E}[\|\widehat{\mathbf{x}}^i_k - \widehat{\mathbf{z}}_k\|] 
\leq\frac{\sum_{\ell=0}^{k}\alpha^2_\ell}{\sum_{\ell=0}^{k}\alpha_\ell}\bigg[N\varsigma\Gamma+\Phi_1\bigg(1+\frac{\Gamma\gamma}{1-\gamma}\bigg)\bigg]\\
&\quad+\frac{1}{\sum_{\ell=0}^{k}\alpha_\ell}\bigg[B_0+\frac{N\varsigma\Gamma\gamma^2}{1-\gamma^2}+\Psi_1\bigg(1+\frac{\Gamma\gamma}{1-\gamma}\bigg)\bigg],
\end{align*}
where $B_0 = \max_i\alpha_0\|\mathbf{x}_0^i-\bar{\mathbf{z}}_0\|$, $\Phi_1>0, \Psi_1>0$ are constants defined in Lemma~\ref{lemma:bound_on_g}, and $\widehat{\mathbf{z}}_k$ is defined in \eqref{eq:def_zhat}.

(2) for any $\mathbf{z}^\star\in\mathcal{X}^\star$, the following relation holds
\begin{align*}
&\mathbb{E}[\|\bar{\mathbf{z}}_{k+1}-\mathbf{z}^\star\|^2|\mathcal{F}_k]\leq\|\bar{\mathbf{z}}_k-\mathbf{z}^\star\|^2\\
&\quad -\frac{2\alpha_k}N(f(\bar{\mathbf{z}}_k) - f^\star) + Z_k,
\end{align*}
where 
\begin{align*}
&Z_k = 2\alpha_k\beta_{1,k}\hat{D}\sqrt{n+2}+4N\varsigma (3\mathcal{K}_1+\bar{\beta}\mathcal{K}_2)\Gamma\gamma^k\alpha_k\\
&+2\mathcal{K}_2\bigg(\frac{2N(2N+\epsilon)\varsigma\Gamma\gamma}{1-\gamma}+\max_{i\in\mathcal{V}}\|\mathbf{x}^i_0-\mathbf{z}^\star\|+{B_1}\bigg)\alpha_k\tilde{\beta}_k\\
&\quad+ 2\mathcal{K}_2\bigg(\frac{(2N+\epsilon)\Gamma\gamma}{1-\gamma} + 2N\bigg)\alpha_k\tilde{\beta}_k\sum_{r=1}^{k-1}\bm{G}_{r-1}\\
&\quad+ \frac{2\mathcal{K}_2}N\alpha_k\tilde{\beta}_k\sum_{r=0}^{k-1}\alpha_r\sum_{i=1}^N\|\mathbf{g}^i(\mathbf{x}^i_r)\|\\
&\quad+2(3\mathcal{K}_1+\bar{\beta}\mathcal{K}_2)\Gamma\alpha_k\sum_{r=1}^k\gamma^{k-r+1}\bm{G}_{r-1}\\
&\quad+2(2\mathcal{K}_1+\bar{\beta}\mathcal{K}_2)\alpha_k\bm{G}_{k-1}+4\varsigma\Gamma\gamma^{k+1}\mathbb{E}[\bm{G}_k|\mathcal{F}_k]\\
&\quad+\frac{2\Gamma}N\sum_{r=1}^k\gamma^{k-r+1}\mathbb{E}[\bm{G}_k|\mathcal{F}_k]\bm{G}_{r-1}\\
&\quad+\frac{4\alpha_k}N\sum_{i=1}^N\mathbb{E}[\|\mathbf{g}^i(\mathbf{x}^i_k)\|\bm{G}_k|\mathcal{F}_k]+\frac5N\mathbb{E}[\bm{G}^2_k|\mathcal{F}_k],
\end{align*}
${B_1}=\max_{i\in\mathcal{V}}\epsilon\|\mathbf{y}^i_0\|+2\sum_{i=1}^N\|\mathbf{x}^i_0-\bar{\mathbf{z}}_0\|$, and $\bar{\beta}$ is the upper bound of $\tilde{\beta}_k$.
\end{Proposition}

\begin{Proof}
For part (1), by the definitions of $\widehat{\mathbf{x}}^i_k$ and $\widehat{\mathbf{z}}_k$, we know that
$\|\widehat{\mathbf{x}}^i_k - \widehat{\mathbf{z}}_k\| \leq \frac{\sum_{\ell=0}^k \alpha_\ell\|\mathbf{x}_\ell^i-\bar{\mathbf{z}}_\ell\|}{\sum_{\ell=0}^k \alpha_\ell}$.
Taking the total expectation and applying Lemma~\ref{lemma:bound_on_consensus}-1), we obtain that for $k\geq1$
\begin{align*}
&\mathbb{E}[\|\widehat{\mathbf{x}}^i_k - \widehat{\mathbf{z}}_k\|] 
\leq \frac{B_0+\sum_{\ell=1}^k \alpha_\ell\mathbb{E}[\|\mathbf{x}_\ell^i-\bar{\mathbf{z}}_\ell\|]}{\sum_{\ell=0}^{k}\alpha_\ell}\\
&\quad\leq\frac{1}{\sum_{\ell=0}^{k}\alpha_\ell}\bigg(B_0+2N\varsigma\Gamma\sum_{\ell=1}^{k}\gamma^\ell\alpha_\ell\\
&\quad\quad+\Gamma\sum_{\ell=1}^{k}\alpha_\ell\sum_{r=1}^{\ell-1}\gamma^{\ell-r}\mathbb{E}[\bm{G}_{r-1}]+ \sum_{\ell=1}^{k}\alpha_\ell\mathbb{E}[\bm{G}_{\ell-1}]\bigg),
\end{align*}
where $B_0 = \max_i\alpha_0\|\mathbf{x}_0^i-\bar{\mathbf{z}}_0\|$ is bounded. 
Since $\alpha_k$ is non-increasing, it follows from Lemma~\ref{lemma:bound_on_g} that
\begin{align}
&\sum_{\ell=1}^k \alpha_\ell\gamma^\ell \leq \frac12\sum_{\ell=1}^{k}\alpha^2_\ell+\frac{\gamma^2}{2(1-\gamma^2)}, \nonumber\\
&\Gamma\sum_{\ell=1}^{k}\alpha_\ell\sum_{r=1}^{\ell-1}\gamma^{\ell-r}\mathbb{E}[\bm{G}_{r-1}] \leq \Gamma\sum_{\ell=1}^{k}\sum_{r=1}^{\ell-1}\gamma^{\ell-r}\alpha_{r-1}\mathbb{E}[\bm{G}_{r-1}]\nonumber\\
&\quad\leq\frac{\Gamma\gamma}{1-\gamma}\sum_{\ell=1}^{k}\alpha_\ell\mathbb{E}[\bm{G}_\ell]\leq\frac{\Phi_1\Gamma\gamma}{1-\gamma}\sum_{\ell=1}^{k}\alpha^2_\ell+\frac{\Psi_1\Gamma\gamma}{1-\gamma}, \nonumber\\
&\sum_{\ell=1}^{k}\alpha_\ell\mathbb{E}[\bm{G}_{\ell-1}] \leq \sum_{\ell=0}^{k}\alpha_\ell\mathbb{E}[\bm{G}_\ell]\leq\Phi_1\sum_{\ell=0}^{k}\alpha^2_\ell+\Psi_1. \label{eq:important_sequence_relations}
\end{align}
Substituting \eqref{eq:important_sequence_relations} to the preceding relation completes the proof of part (1).

For part (2), considering \eqref{eq:d-dgd}, and the fact that $A$ is column-stochastic, we have
$\bar{\mathbf{z}}_{k+1} = \bar{\mathbf{z}}_k+ \frac1N\sum_{i=1}^Ng^i_k$.
Then, for any $\mathbf{z}^\star\in\mathcal{X}^\star$, it follows that
\begin{subequations}\label{eq:last_term}
\begin{align}
&\|\bar{\mathbf{z}}_{k+1}-\mathbf{z}^\star\|^2\leq\|\bar{\mathbf{z}}_k-\mathbf{z}^\star\|^2 + \frac{\bm{G}^2_k}{N^2}+\frac2N\sum_{i=1}^N\langle g^i_k,\bar{\mathbf{z}}_k-\mathbf{z}^\star\rangle\nonumber\\
&\quad=\|\bar{\mathbf{z}}_k-\mathbf{z}^\star\|^2 + \frac{\bm{G}^2_k}{N^2}\label{eq:term1}\\
&\quad\quad+\frac2N\sum_{i=1}^N(g^i_k+\alpha_k\mathbf{g}^i(\mathbf{x}^i_k)),\bar{\mathbf{z}}_k-\mathbf{z}^\star\rangle\label{eq:term3}\\
&\quad\quad-\frac{2\alpha_k}N\sum_{i=1}^N\langle \mathbf{g}^i(\mathbf{x}^i_k), \bar{\mathbf{z}}_k-\mathbf{z}^\star\rangle.\label{eq:term2}
\end{align}
\end{subequations}
For the second term in \eqref{eq:term1}, we have that $\mathbb{E}[\frac{\bm{G}^2_k}{N^2}|\mathcal{F}_k]\leq \frac1N\mathbb{E}[\bm{G}^2_k|\mathcal{F}_k]$.

For term \eqref{eq:term3}, it can be expanded as
\begin{subequations}\label{eq:last_term_break}
\begin{align}
&\sum_{i=1}^N\langle g^i_k+\alpha_k\mathbf{g}^i(\mathbf{x}^i_k),\bar{\mathbf{z}}_k-\mathbf{z}^\star\rangle\nonumber\\
=&\sum_{i=1}^N\langle g^i_k+\alpha_k\mathbf{g}^i(\mathbf{x}^i_k),\bar{\mathbf{z}}_k-\bar{\mathbf{z}}_{k+1}\rangle\label{eq:last_term_break_one}\\
&+\sum_{i=1}^N\langle g^i_k+\alpha_k\mathbf{g}^i(\mathbf{x}^i_k),\bar{\mathbf{z}}_{k+1}-\mathbf{x}^i_{k+1}\rangle\label{eq:last_term_break_two}\\
&+\sum_{i=1}^N\langle g^i_k+\alpha_k\mathbf{g}^i(\mathbf{x}^i_k),\mathbf{x}^i_{k+1}-\mathbf{z}^\star\rangle.\label{eq:last_term_break_three}
\end{align}
\end{subequations}
For \eqref{eq:last_term_break_one}, we have
\begin{align}
&\sum_{i=1}^N\mathbb{E}[\langle g^i_k+\alpha_k\mathbf{g}^i(\mathbf{x}^i_k),\bar{\mathbf{z}}_k-\bar{\mathbf{z}}_{k+1}\rangle|\mathcal{F}_k]\nonumber\\
&\quad\leq\frac1N\mathbb{E}[\bm{G}^2_k|\mathcal{F}_k]+\frac{\alpha_k}N\sum_{i=1}^N\mathbb{E}[\|\mathbf{g}^i(\mathbf{x}^i_k)\|\bm{G}_k|\mathcal{F}_k]\nonumber\\
&\quad\leq \mathbb{E}[\bm{G}^2_k|\mathcal{F}_k]+\alpha_k\sum_{i=1}^N\mathbb{E}[\|\mathbf{g}^i(\mathbf{x}^i_k)\|\bm{G}_k|\mathcal{F}_k].\label{eq:last_term_break_1}
\end{align}
For \eqref{eq:last_term_break_two}, we have
$\sum_{i=1}^N\langle g^i_k+\alpha_k\mathbf{g}^i(\mathbf{x}^i_k),\bar{\mathbf{z}}_{k+1}-\mathbf{x}^i_{k+1}\rangle
\leq(\bm{G}_k+\alpha_k\sum_{i=1}^N\|\mathbf{g}^i(\mathbf{x}^i_k)\|)(2N\varsigma\Gamma\gamma^{k+1}
+\Gamma\sum_{r=1}^k\gamma^{k-r+1}\bm{G}_{r-1}+ \bm{G}_k)$,
where Lemma~\ref{lemma:bound_on_consensus}-(1) was substituted. Hence, we obtain
\begin{align}
&\sum_{i=1}^N\mathbb{E}[\langle g^i_k+\alpha_k\mathbf{g}^i(\mathbf{x}^i_k),\bar{\mathbf{z}}_{k+1}-\mathbf{x}^i_{k+1}\rangle|\mathcal{F}_k]\nonumber\\
&\quad\leq2N^2\varsigma\mathcal{K}_1\Gamma\gamma^{k+1}\alpha_k+2N\varsigma\Gamma\gamma^{k+1}\mathbb{E}[\bm{G}_k|\mathcal{F}_k]\nonumber\\
&\quad\quad+N\mathcal{K}_1\Gamma\alpha_k\sum_{r=1}^k\gamma^{k-r+1}\bm{G}_{r-1}\nonumber\\
&\quad\quad+\Gamma\sum_{r=1}^k\gamma^{k-r+1}\mathbb{E}[\bm{G}_k|\mathcal{F}_k]\bm{G}_{r-1}\nonumber\\
&\quad\quad+\alpha_k\sum_{i=1}^N\mathbb{E}[\|\mathbf{g}^i(\mathbf{x}^i_k)\|\bm{G}_k|\mathcal{F}_k]+\mathbb{E}[\bm{G}^2_k|\mathcal{F}_k].\label{eq:last_term_break_2}
\end{align}
For \eqref{eq:last_term_break_three}, it follows from \cite[Lemma 1-(a)]{Nedic2010} that
\begin{align}\label{eq:last_term_break_3}
\langle g^i_k+\alpha_k\mathbf{g}^i(\mathbf{x}^i_k),\mathbf{x}^i_{k+1}-\mathbf{z}^\star\rangle \leq 0.
\end{align}
Thus, taking the conditional expectation on $\mathcal{F}_k$ in \eqref{eq:last_term_break} and substituting \eqref{eq:last_term_break_1}, \eqref{eq:last_term_break_2} and \eqref{eq:last_term_break_3}, we obtain
\begin{align}
&\sum_{i=1}^N\mathbb{E}[\langle g^i_k+\alpha_k\mathbf{g}^i(\mathbf{x}^i_k),\bar{\mathbf{z}}_k-\mathbf{z}^\star\rangle|\mathcal{F}_k]\nonumber\\
&\quad\leq2N^2\varsigma\mathcal{K}_1\Gamma\gamma^{k+1}\alpha_k+2N\varsigma\Gamma\gamma^{k+1}\mathbb{E}[\bm{G}_k|\mathcal{F}_k]\nonumber\\
&\quad\quad+N\mathcal{K}_1\Gamma\alpha_k\sum_{r=1}^k\gamma^{k-r+1}\bm{G}_{r-1}\nonumber\\
&\quad\quad+\Gamma\sum_{r=1}^k\gamma^{k-r+1}\mathbb{E}[\bm{G}_k|\mathcal{F}_k]\bm{G}_{r-1}\nonumber\\
&\quad\quad+2\alpha_k\sum_{i=1}^N\mathbb{E}[\|\mathbf{g}^i(\mathbf{x}^i_k)\|\bm{G}_k|\mathcal{F}_k]+2\mathbb{E}[\bm{G}^2_k|\mathcal{F}_k].\label{eq:last_term_combine}
\end{align}
For \eqref{eq:term2}, from Lemma~\ref{lemma:property_f_mu}-(2),
$\sum_{i=1}^N\mathbb{E}[\langle \mathbf{g}^i(\mathbf{x}^i_k), \bar{\mathbf{z}}_k-\mathbf{z}^\star\rangle|\mathcal{F}_k]
=\sum_{i=1}^N\langle\nabla f_{i,\beta_{1,k}}(\mathbf{x}^i_k)+\tilde{\beta}_k\hat{D}\mathbf{v},\bar{\mathbf{z}}_k-\mathbf{z}^\star\rangle$.
Denote $\hat{D}\|\mathbf{v}\|$ by $\mathcal{K}_2$, we have
\begin{align}
&\langle\nabla f_{i,\beta_{1,k}}(\mathbf{x}^i_k)+\tilde{\beta}_k\hat{D}\mathbf{v},\bar{\mathbf{z}}_k-\mathbf{z}^\star\rangle\nonumber\\
&= \langle\nabla f_{i,\beta_{1,k}}(\mathbf{x}^i_k)+\tilde{\beta}_k\hat{D}\mathbf{v},\bar{\mathbf{z}}_k-\mathbf{x}^i_k\rangle\nonumber\\
&\quad+ \langle\nabla f_{i,\beta_{1,k}}(\mathbf{x}^i_k)+\tilde{\beta}_k\hat{D}\mathbf{v},\mathbf{x}^i_k-\mathbf{z}^\star\rangle\nonumber\\
&\geq -\|\nabla f_{i,\beta_{1,k}}(\mathbf{x}^i_k)\|\|\mathbf{x}^i_k-\bar{\mathbf{z}}_k\|-\tilde{\beta}_k\mathcal{K}_2\|\mathbf{x}^i_k-\bar{\mathbf{z}}_k\|\nonumber\\
&\quad+ f_{i,\beta_{1,k}}(\mathbf{x}^i_k) - f_{i,\beta_{1,k}}(\mathbf{z}^\star)-\tilde{\beta}_k\mathcal{K}_2\|\mathbf{x}^i_k-\mathbf{z}^\star\|\nonumber\\
&\geq f_{i,\beta_{1,k}}(\bar{\mathbf{z}}_k) - f_{i,\beta_{1,k}}(\mathbf{z}^\star)-\tilde{\beta}_k\mathcal{K}_2\|\mathbf{x}^i_k-\mathbf{z}^\star\|\nonumber\\
&\quad-(2\mathcal{K}_1+\tilde{\beta}_k\mathcal{K}_2)\|\mathbf{x}^i_k-\bar{\mathbf{z}}_k\|\nonumber\\
&\geq f_i(\bar{\mathbf{z}}_k) - f_i(\mathbf{z}^\star) - \beta_{1,k}\hat{D}\sqrt{n+2}-\tilde{\beta}_k\mathcal{K}_2\|\mathbf{x}^i_k\nonumber\\
&\quad-\mathbf{z}^\star\|-(2\mathcal{K}_1+\tilde{\beta}_k\mathcal{K}_2)\|\mathbf{x}^i_k-\bar{\mathbf{z}}_k\|.\label{eq:g_times_zbar_minus_z_star}
\end{align}
Considering the term $\|\mathbf{x}^i_k-\mathbf{z}^\star\|$, it follows that
$\|\mathbf{x}^i_k-\mathbf{z}^\star\|
\leq\sum_{j=1}^N [A_r]_{ij}\|\mathbf{x}^j_{k-1}-\mathbf{z}^\star\| + \epsilon \|\mathbf{y}^i_{k-1}\| + \alpha_{k-1} \|\mathbf{g}^i(\mathbf{x}^i_{k-1})\|
\leq\sum_{j=1}^N [A_r]_{ij}\|\mathbf{x}^i_{k-1}-\mathbf{z}^\star\| + \epsilon \|\mathbf{y}^i_{k-1}\| + \alpha_{k-1} \|\mathbf{g}^i(\mathbf{x}^i_{k-1})\|
+\sum_{j=1}^N [A_r]_{ij}\|\mathbf{x}^i_{k-1}-\mathbf{x}^j_{k-1}\|
\leq\|\mathbf{x}^i_{k-1}-\mathbf{z}^\star\| + \epsilon \|\mathbf{y}^i_{k-1}\| + \alpha_{k-1} \|\mathbf{g}^i(\mathbf{x}^i_{k-1})\|
+2\sum_{i=1}^N\|\mathbf{x}^i_{k-1}-\bar{\mathbf{z}}_{k-1}\|$.

Applying the above relation recursively yields
$\|\mathbf{x}^i_k-\mathbf{z}^\star\| \leq \|\mathbf{x}^i_0-\mathbf{z}^\star\|+\epsilon \sum_{\tau=0}^{k-1}\|\mathbf{y}^i_\tau\| + \sum_{\tau=0}^{k-1}\alpha_\tau \|\mathbf{g}^i(\mathbf{x}^i_\tau)\|
+2\sum_{\tau=0}^{k-1}\sum_{i=1}^N \|\mathbf{x}^i_\tau-\bar{\mathbf{z}}_\tau\|$.

Thus, substituting the above result to \eqref{eq:g_times_zbar_minus_z_star} gives
$\langle\nabla f_{i,\beta_{1,k}}(\mathbf{x}^i_k)+\tilde{\beta}_k\hat{D}\mathbf{v},\bar{\mathbf{z}}_k-\mathbf{z}^\star\rangle\geq f_i(\bar{\mathbf{z}}_k) - f_i(\mathbf{z}^\star)
-\tilde{\beta}_k\mathcal{K}_2(\|\mathbf{x}^i_0-\mathbf{z}^\star\|+\epsilon \sum_{\tau=0}^{k-1}\|\mathbf{y}^i_\tau\| + \sum_{\tau=0}^{k-1}\alpha_\tau \|\mathbf{g}^i(\mathbf{x}^i_\tau)\|
+2\sum_{\tau=0}^{k-1}\sum_{i=1}^N \|\mathbf{x}^i_\tau-\bar{\mathbf{z}}_\tau\|)
-(2\mathcal{K}_1+\tilde{\beta}_k\mathcal{K}_2)\|\mathbf{x}^i_k-\bar{\mathbf{z}}_k\| - \beta_{1,k}\hat{D}\sqrt{n+2}$.
Applying Lemma~\ref{lemma:bound_on_consensus} and noting that $\tilde{\beta}_k$ is bounded (where its upper bound is denoted by $\bar{\beta}$), we obtain that
\begin{align}
&\frac{2\alpha_k}N\sum_{i=1}^N\mathbb{E}[\langle \mathbf{g}^i(\mathbf{x}^i_k), \bar{\mathbf{z}}_k-\mathbf{z}^\star\rangle|\mathcal{F}_k]\geq \frac{2\alpha_k}N(f(\bar{\mathbf{z}}_k) - f^\star)\nonumber\\
&\quad- 2\alpha_k\beta_{1,k}\hat{D}\sqrt{n+2}-2\mathcal{K}_2\bigg(\frac{2N(2N+\epsilon)\varsigma\Gamma\gamma}{1-\gamma}\nonumber\\
&\quad\quad+\max_{i\in\mathcal{V}}\|\mathbf{x}^i_0-\mathbf{z}^\star\|+{B_1}\bigg)\alpha_k\tilde{\beta}_k\nonumber\\
&\quad- 2\mathcal{K}_2\bigg(\frac{(2N+\epsilon)\Gamma\gamma}{1-\gamma} + 2N\bigg)\alpha_k\tilde{\beta}_k\sum_{r=1}^{k-1}\bm{G}_{r-1}\nonumber\\
&\quad- \frac{2\mathcal{K}_2}N\alpha_k\tilde{\beta}_k\sum_{r=0}^{k-1}\alpha_r\sum_{i=1}^N\|\mathbf{g}^i(\mathbf{x}^i_r)\|\nonumber\\
&\quad-4N\varsigma (2\mathcal{K}_1+\bar{\beta}\mathcal{K}_2)\Gamma\gamma^k\alpha_k-2(2\mathcal{K}_1+\bar{\beta}\mathcal{K}_2)\alpha_k\bm{G}_{k-1}\nonumber\\
&\quad-2(2\mathcal{K}_1+\bar{\beta}\mathcal{K}_2)\Gamma\alpha_k\sum_{r=1}^{k-1}\gamma^{k-r}\bm{G}_{r-1},\label{eq:E_g_z_minus_xstar_2}
\end{align}
where ${B_1}=\max_{i\in\mathcal{V}}\epsilon\|\mathbf{y}^i_0\|+2\sum_{i=1}^N\|\mathbf{x}^i_0-\bar{\mathbf{z}}_0\|$.

Taking the conditional expectation on $\mathcal{F}_k$ in \eqref{eq:last_term}, and substituting \eqref{eq:last_term_combine} and \eqref{eq:E_g_z_minus_xstar_2} gives the result of part (2).
\end{Proof}

\subsection{Main Results}\label{subsec:main_results}

In this subsection, we present the main convergence results of our proposed algorithm, including convergence under non-summable and square-summable step-size condition (Theorem~\ref{theorem:non_summable_square_summable}), convergence under non-summable step-size condition only (Theorem~\ref{theorem:non_summable}), and the convergence rate analysis (Corollary~\ref{corollary:convergence_rate}).

Our first result demonstrates the standard convergence results under non-summable and square-summable step-size condition.

\begin{Theorem}\label{theorem:non_summable_square_summable}
Suppose Assumptions~\ref{assumption_graph}, \ref{assumption_local_f} and \ref{assumption_random_variables} hold. Let $\{\widehat{\mathbf{x}}^i_k\}_{k\geq0}$ be the sequence generated by \eqref{eq:algorithm} with a non-increasing step-size sequence $\{\alpha_k\}_{k\geq0}$ satisfying 
\begin{align*}
\sum_{k=0}^\infty\alpha_k= \infty, \quad\sum_{k=0}^\infty \alpha^2_k< \infty.
\end{align*}
Let $\epsilon$ be the constant such that $0<\epsilon <\min(\bar{\epsilon},\frac{1-\gamma}{2\sqrt{3}N\Gamma\gamma})$, where $\Gamma$ and $\gamma$ are some constants, and $\bar{\epsilon} = (\frac{1-|\lambda_3|}{20+8N})^N$ with $\lambda_3$ being the third largest eigenvalue of the weighting matrix $A$ in \eqref{eq:d-dgd} by setting $\epsilon =0$. Let $\beta_{1,k}$ and $\tilde{\beta}_k$ defined in \eqref{eq:ratio_two_sequence} satisfy $\sum_{k=0}^\infty\beta_{1,k}\alpha_k < \infty$ and $\sum_{k=0}^\infty\tilde{\beta}_k < \infty$. Then, for $\forall i\in\mathcal{V}$, we have $\{\widehat{\mathbf{x}}^i_k\}_{k\geq0}$ converges a.s. to an optimizer $\mathbf{x}^\star\in\mathcal{X}^\star$\footnote{In this paper, `a.s.' is meant for `almost surely'. For a sequence of random vectors $\{a_k\}_{k\geq0}$, we say that $a_k$ converges to $a$ almost surely, if $\mathbb{P}(\lim_{k\to\infty}a_k=a)=1$, \textit{i.e}, the probability of $\lim_{k\to\infty}a_k=a$ is 1.}.
\end{Theorem}

\begin{Proof}
we proceed to the proof by showing (A) the convergence of $\widehat{\mathbf{x}}^i_k$ to $\mathbf{x}^i_k$, (B) the convergence of $\mathbf{x}^i_k$ to $\bar{\mathbf{z}}_k$, and (C) the convergence of $\bar{\mathbf{z}}_k$ to an optimizer $\mathbf{x}^\star\in\mathcal{X}^\star$ under appropriate conditions.

(A) Convergence of $\widehat{\mathbf{x}}^i_k$ to $\mathbf{x}^i_k$:

Suppose $\{\mathbf{x}^i_k\}_{k\geq0}$ converges a.s. to some point $\tilde{\mathbf{x}}$, \textit{i.e.}, $\mathbb{P}(\lim_{k\to\infty}\mathbf{x}^i_k = \tilde{\mathbf{x}})=1$. By definition of $\widehat{\mathbf{x}}^i_k=\frac{\sum_{\ell=0}^k\alpha_\ell\mathbf{x}^i_\ell}{\sum_{\ell=0}^k\alpha_\ell}$ in \eqref{eq:update_xhat}, it follows from Lemma~\ref{theorem:Stolz-Cesaro} with $a_k = \alpha_k\mathbf{x}^i_k$ and $b_k = \alpha_k$ that $\mathbb{P}(\lim_{k\to\infty}\widehat{\mathbf{x}}^i_k=\lim_{k\to\infty}\mathbf{x}^i_k = \tilde{\mathbf{x}})=1$. Hence, we obtain $\{\widehat{\mathbf{x}}^i_k\}_{k\geq0}$ converges a.s. to the same point $\tilde{\mathbf{x}}$.

(B) Convergence of $\mathbf{x}^i_k$ to $\bar{\mathbf{z}}_k$:

Squaring both sides of Lemma~\ref{lemma:bound_on_consensus}-(1), taking the total expectation, and summing over from $k=1$ to infinity, we obtain
\begin{align*}
\sum_{k=1}^\infty\mathbb{E}[\|\mathbf{x}^i_k-\bar{\mathbf{z}}_k\|^2]&\leq \frac{12N^2\varsigma^2\Gamma^2\gamma^2}{1-\gamma^2}+\bigg(\frac{3\Psi_2\Gamma^2\gamma^2}{(1-\gamma)^2}+ 3\Psi_2\bigg)\\
&\quad+\bigg(\frac{3\Phi_2\Gamma^2\gamma^2}{(1-\gamma)^2}+ 3\Phi_2\bigg)\sum_{k=0}^\infty\alpha^2_k,
\end{align*}
where $(\sum_{r=1}^{k-1}\gamma^{k-r}\bm{G}_{r-1})^2\leq\frac{\gamma}{1-\gamma}\sum_{r=1}^{k-1}\gamma^{k-r}\bm{G}^2_{r-1}$, $\sum_{k=1}^\infty\sum_{r=1}^{k-1}\gamma^{k-r}\mathbb{E}[\bm{G}^2_{r-1}] \leq \frac{\gamma}{1-\gamma}\sum_{k=0}^\infty\mathbb{E}[\bm{G}_k^2]$ and Lemma~\ref{lemma:bound_on_g} have been applied. 
As the step-size is square-summable, we obtain $\sum_{k=1}^\infty\mathbb{E}[\|\mathbf{x}^i_k-\bar{\mathbf{z}}_k\|^2]< \infty$. By the monotone convergence theorem, it follows that $\mathbb{E}[\sum_{k=1}^\infty\|\mathbf{x}^i_k-\bar{\mathbf{z}}_k\|^2]=\sum_{k=1}^\infty\mathbb{E}[\|\mathbf{x}^i_k-\bar{\mathbf{z}}_k\|^2]< \infty$, which implies $\{\mathbf{x}^i_k-\bar{\mathbf{z}}_k\}_{k\geq0}$ converges a.s. to $0$.

(C) Convergence of $\bar{\mathbf{z}}_k$ to an optimizer $\mathbf{x}^\star\in\mathcal{X}^\star$:

Finally, we will show that $\{\bar{\mathbf{z}}_k\}_{k\geq0}$ indeed has a limit, and converges to an optimizer $\mathbf{x}^\star\in\mathcal{X}^\star$. The proof of this part is based on the Robbins-Siegmund's Lemma\cite{Polyak1987} as quoted below for completeness.

\begin{Lemma}\label{lemma:Robbins-Siegmund}
(Robbins-Siegmund's Lemma) Let $u_k, v_k, w_k,\eta_k$ be non-negative random variables satisfying that
\begin{gather*}
\mathbb{E}[u_{k+1}|\mathcal{F}_k] \leq (1+\eta_k)u_k - v_k + w_k \quad \text{a.s.},\\
\sum_{k=0}^\infty \eta_k < \infty \quad \text{a.s.}, \qquad \sum_{k=0}^\infty w_k < \infty \quad \text{a.s.},
\end{gather*}
where $\mathbb{E}[u_{k+1}|\mathcal{F}_k]$ is the conditional expectation for the given $u_0,\ldots,u_k,v_0,\ldots,v_k,w_0,\ldots,w_k,\eta_0,\ldots,\eta_k$. Then
\begin{enumerate}
\item $\begin{aligned}[t]\{u_k\}_{k\geq0}\end{aligned}$ converges a.s.;
\item $\begin{aligned}[t]\sum_{k=0}^\infty v_k < \infty\end{aligned}$ a.s.
\end{enumerate}
\end{Lemma}

From Proposition~\ref{prop:consensus}-(2), we have that for any $\mathbf{z}^\star\in\mathcal{X}^\star$,
$\mathbb{E}[\|\bar{\mathbf{z}}_{k+1}-\mathbf{z}^\star\|^2|\mathcal{F}_k]\leq\|\bar{\mathbf{z}}_k-\mathbf{z}^\star\|^2 
-\frac{2\alpha_k}N(f(\bar{\mathbf{z}}_k) - f^\star) + Z_k$.
To invoke Lemma~\ref{lemma:Robbins-Siegmund}, it suffices to show that $\sum_{k=0}^\infty Z_k<\infty$, a.s.

Now, taking the total expectation for $Z_k$ and summing over from $k=1$ to infinity, we have
\begin{align*}
&\sum_{k=1}^\infty\mathbb{E}[Z_k] \leq 2\hat{D}\sqrt{n+2}\sum_{k=1}^\infty\alpha_k\beta_{1,k}+2\alpha_0\mathcal{K}_2\bigg({B_1}\\
&\quad\quad+\frac{2N(2N+\epsilon)\varsigma\Gamma\gamma}{1-\gamma}+\max_{i\in\mathcal{V}}\|\mathbf{x}^i_0-\mathbf{z}^\star\|\bigg)\sum_{k=0}^\infty\tilde{\beta}_k\\
&\quad+ 2\mathcal{K}_2\bigg[\mathcal{K}_1+ \Phi_1\bigg(\frac{(2N+\epsilon)\Gamma\gamma}{1-\gamma} + 2N\bigg)\bigg]\sum_{k=0}^\infty\tilde{\beta}_k\sum_{k=0}^\infty\alpha^2_k\\
&\quad+2\Phi_1(2\mathcal{K}_1+\bar{\beta}\mathcal{K}_2)\sum_{k=0}^\infty\alpha^2_k+\frac{5\Phi_2+4\Phi_3}N\sum_{k=1}^\infty\alpha_k^2\\
&\quad+2\varsigma (\Phi_2+3N\mathcal{K}_1+N\bar{\beta}\mathcal{K}_2)\Gamma\sum_{k=1}^\infty\alpha^2_k\\
&\quad+2\bigg(\Phi_1(3\mathcal{K}_1+\bar{\beta}\mathcal{K}_2)+\frac{\Phi_2}N\bigg)\frac{\Gamma\gamma}{1-\gamma}\sum_{k=1}^\infty\alpha^2_k\\
&\quad+2\mathcal{K}_2\Psi_1\bigg(\frac{(2N+\epsilon)\Gamma\gamma}{1-\gamma} + 2N\bigg)\sum_{k=0}^\infty\tilde{\beta}_k+2\Psi_1(2\mathcal{K}_1\\
&\quad\quad+\bar{\beta}\mathcal{K}_2)+2\varsigma\Gamma\bigg(\frac{(1+3N\mathcal{K}_1+N\bar{\beta}\mathcal{K}_2)\gamma^2}{1-\gamma^2}+\Psi_2\bigg)\\
&\quad+2\bigg(\Psi_1(3\mathcal{K}_1+\bar{\beta}\mathcal{K}_2)+\frac{\Psi_2}N\bigg)\frac{\Gamma\gamma}{1-\gamma}+\frac{5\Psi_2+4\Psi_3}N,
\end{align*}
where we applied $\mathbb{E}[\mathbb{E}[\bm{G}_k|\mathcal{F}_k]\bm{G}_{r-1}]\leq\sqrt{\mathbb{E}[\bm{G}^2_k]\mathbb{E}[\bm{G}^2_{r-1}]}\leq\frac12(\mathbb{E}[\bm{G}^2_k]+\mathbb{E}[\bm{G}^2_{r-1}])$ based on Cauchy-Schwarz inequality, the results in \eqref{eq:important_sequence_relations}, and $\sum_{k=0}^\infty\alpha_k\tilde{\beta}_k\sum_{r=0}^{k-1}\mathbb{E}[\bm{G}_r] \leq \sum_{k=0}^\infty\tilde{\beta}_k\sum_{r=0}^{k-1}\alpha_r\mathbb{E}[\bm{G}_r]\leq\sum_{k=0}^\infty\tilde{\beta}_k\sum_{k=0}^\infty\alpha_k\mathbb{E}[\bm{G}_k]\leq\Phi_1\sum_{k=0}^\infty\tilde{\beta}_k\sum_{k=0}^\infty\alpha_k^2+\Psi_1\sum_{k=0}^\infty\tilde{\beta}_k$.
Since $\sum_{k=0}^\infty\beta_{1,k}\alpha_k < \infty$, $\sum_{k=0}^\infty\tilde{\beta}_k < \infty$ and $\sum_{k=0}^\infty \alpha^2_k< \infty$, by the monotone convergence theorem, we have $\mathbb{E}[\sum_{k=1}^\infty Z_k]=\sum_{k=1}^\infty \mathbb{E}[Z_k]<\infty$, which proves that $\sum_{k=1}^\infty Z_k<\infty$ a.s. 

Invoking Lemma~\ref{lemma:Robbins-Siegmund}, we obtain that 
\begin{subequations}
\begin{align}
\forall \mathbf{z}^\star\in\mathcal{X}^\star, \{\|\bar{\mathbf{z}}_k-\mathbf{z}^\star\|^2\}_{k\geq0} \text{ converges a.s. } \label{eq:result1}\\
\sum_{k=0}^\infty\alpha_k(f(\bar{\mathbf{z}}_k) - f^\star)<\infty \text{ a.s.} \label{eq:result2}
\end{align}
\end{subequations}
Since $f(\bar{\mathbf{z}}_k) - f^\star\geq0$, and the step-size is non-summable, it follows from \eqref{eq:result2} that $\liminf_{k\to\infty}f(\bar{\mathbf{z}}_k) = f^\star$ a.s.
Let $\{\bar{\mathbf{z}}_{k_1}\}_{k_1\geq0}$ be a subsequence of $\{\bar{\mathbf{z}}_k\}_{k\geq0}$ such that
\begin{align} 
\lim_{k_1\to\infty}f(\bar{\mathbf{z}}_{k_1}) = \liminf_{k\to\infty}f(\bar{\mathbf{z}}_k) = f^\star \text{ a.s.} \label{eq:limit_subsequence}
\end{align} 
From \eqref{eq:result1}, the sequence $\{\bar{\mathbf{z}}_k\}_{k\geq0}$ is bounded a.s.
Without loss of generality, we may assume that $\{\bar{\mathbf{z}}_{k_1}\}_{k_1\geq0}$ converges a.s. to some $\mathbf{x}^\star$ (if not, we may choose one such subsequence). Due to the continuity of $f$, we have $f(\bar{\mathbf{z}}_{k_1})$ converges to $f(\mathbf{x}^\star)$ a.s., which by \eqref{eq:limit_subsequence} implies that $f(\mathbf{x}^\star)=f^\star$, \textit{i.e.}, $\mathbf{x}^\star\in\mathcal{X}^\star$. Then we let $\mathbf{z}^\star = \mathbf{x}^\star$ in \eqref{eq:result1} and consider the sequence $\{\|\bar{\mathbf{z}}_k-\mathbf{x}^\star\|^2\}_{k\geq0}$. It converges a.s., and its subsequence $\{\|\bar{\mathbf{z}}_{k_1}-\mathbf{x}^\star\|^2\}_{k\geq0}$ converges a.s. to $0$. Thus, we have $\{\bar{\mathbf{z}}_k\}_{k\geq0}$ converges a.s. to $\mathbf{x}^\star$.

Therefore, combining the arguments of (A), (B) and (C), we complete the proof of Theorem~\ref{theorem:non_summable_square_summable}.
\end{Proof}


Our second result removes the square-summable step-size condition, and shows the convergence of $\mathbb{E}[f(\widehat{\mathbf{x}}^i_k)]$ to the optimal value.

\begin{Theorem}\label{theorem:non_summable}
Suppose Assumptions~\ref{assumption_graph}, \ref{assumption_local_f} and \ref{assumption_random_variables} hold. Let $\{\widehat{\mathbf{x}}^i_k\}_{k\geq0}$ be the sequence generated by \eqref{eq:algorithm} with a non-increasing step-size sequence $\{\alpha_k\}_{k\geq0}$ satisfying 
\begin{align*}
\sum_{k=0}^\infty\alpha_k= \infty, \quad\lim_{k\to\infty} \alpha_k=\alpha_\infty.
\end{align*} 
Let $\epsilon$ be the constant such that $0<\epsilon < \min(\bar{\epsilon},\frac{1-\gamma}{2\sqrt{3}N\Gamma\gamma})$, where $\Gamma$ and $\gamma$ are some constants, and $\bar{\epsilon} = (\frac{1-|\lambda_3|}{20+8N})^N$ with $\lambda_3$ being the third largest eigenvalue of the weighting matrix $A$ in \eqref{eq:d-dgd} by setting $\epsilon =0$. Let $\beta_{1,k}$ and $\tilde{\beta}_k$ defined in \eqref{eq:ratio_two_sequence} satisfy $\lim_{k\to\infty}\beta_{1,k} = 0$ and $\sum_{k=0}^\infty\tilde{\beta}_k < \infty$. Then, for any $\mathbf{z}^\star\in\mathcal{X}^\star$, we have
\begin{align*}
&\limsup_{k\to\infty}\mathbb{E}[f(\widehat{\mathbf{x}}^i_k)] - f^\star
\leq \alpha_\infty\sum_{k=0}^\infty\tilde{\beta}_kN\mathcal{K}_2\bigg[\mathcal{K}_1\\
&\quad+\Phi_1\bigg(\frac{(2N+\epsilon)\Gamma\gamma}{1-\gamma} + 2N\bigg)\bigg]+\alpha_\infty \bigg[2.5\Phi_2+2\Phi_3\\
&\quad+\frac{(N\Phi_1(3\mathcal{K}_1+\bar{\beta}\mathcal{K}_2)+\Phi_2+\hat{D}\Phi_1)\Gamma\gamma}{1-\gamma}+\hat{D}\Phi_1\\
&\quad+N\varsigma (\Phi_2+N(3\mathcal{K}_1+\bar{\beta}\mathcal{K}_2)+\hat{D})\Gamma\bigg],
\end{align*}
where $f^\star$ is the optimal value of the problem, \textit{i.e.}, $f^\star = \min_{\mathbf{z}^\star\in\mathcal{X}^\star}f(\mathbf{z}^\star)$, $\mathcal{K}_1$, $\mathcal{K}_2$, $\Phi_1$, $\Phi_2$, $\Phi_3$, $\Gamma$ and $\gamma$ are positive constants, $\bar{\beta}$ is the upper bound of $\tilde{\beta}_k$, and $\varsigma = \max\{\|\mathbf{x}^i_0\|,\|\mathbf{y}^i_0\|,i\in\mathcal{V}\}$.
\end{Theorem}

\begin{Proof}
Taking the total expection for the result in Proposition~\ref{prop:consensus}-(2) and re-arranging the terms, we have
$\alpha_k(\mathbb{E}[f(\bar{\mathbf{z}}_k)] - f^\star)
\leq \frac{N}2(\mathbb{E}[\|\bar{\mathbf{z}}_k-\mathbf{z}^\star\|^2]-\mathbb{E}[\|\bar{\mathbf{z}}_{k+1}-\mathbf{z}^\star\|^2|]) + \frac{N}2\mathbb{E}[Z_k]$.
Summing over from $k = 0$ to $t-1$,
we have
\begin{align}
&\sum_{k=0}^{t-1}\alpha_k(\mathbb{E}[f(\bar{\mathbf{z}}_k)] - f^\star)\leq \frac{N}2\sum_{k=0}^{t-1}\mathbb{E}[Z_k]+\frac{N}2\|\bar{\mathbf{z}}_0-\mathbf{z}^\star\|^2\nonumber\\
&\leq N\hat{D}\sqrt{n+2}\sum_{k=0}^{t-1}\alpha_k\beta_{1,k}+N\mathcal{K}_2\bigg(\frac{2N(2N+\epsilon)\varsigma\Gamma\gamma}{1-\gamma}\nonumber\\
&\quad\quad+\max_{i\in\mathcal{V}}\|\mathbf{x}^i_0-\mathbf{z}^\star\|+{B_1}\bigg)\sum_{k=0}^{t-1}\alpha_k\tilde{\beta}_k+ N\mathcal{K}_2\bigg[\mathcal{K}_1\nonumber\\
&\quad\quad+ \Phi_1\bigg(\frac{(2N+\epsilon)\Gamma\gamma}{1-\gamma} + 2N\bigg)\bigg]\sum_{k=0}^{t-1}\tilde{\beta}_k\sum_{k=0}^{t-1}\alpha^2_k\nonumber\\
&\quad+\bigg[\frac{(N\Phi_1(3\mathcal{K}_1+\bar{\beta}\mathcal{K}_2)+\Phi_2)\Gamma\gamma}{1-\gamma}+2.5\Phi_2+2\Phi_3\nonumber\\
&\quad\quad+N\varsigma (\Phi_2+N(3\mathcal{K}_1+\bar{\beta}\mathcal{K}_2))\Gamma\bigg]\sum_{k=0}^{t-1}\alpha^2_k\nonumber\\
&\quad+N\mathcal{K}_2\Psi_1\bigg(\frac{(2N+\epsilon)\Gamma\gamma}{1-\gamma} + 2N\bigg)\sum_{k=0}^\infty\tilde{\beta}_k\nonumber\\
&\quad+N\varsigma\Gamma\bigg(\frac{(1+3N\mathcal{K}_1+N\bar{\beta}\mathcal{K}_2)\gamma^2}{1-\gamma^2}+\Psi_2\bigg)\nonumber\\
&\quad+\frac{(N\Psi_1(3\mathcal{K}_1+\bar{\beta}\mathcal{K}_2)+\Psi_2)\Gamma\gamma}{1-\gamma}+\frac{N}2\|\bar{\mathbf{z}}_0-\mathbf{z}^\star\|^2\nonumber\\
&\quad+N\Psi_1(2\mathcal{K}_1+\bar{\beta}\mathcal{K}_2)+2.5\Psi_2+2\Psi_3.\label{eq:E_f-f_mu}
\end{align}

Dividing both sides of \eqref{eq:E_f-f_mu} by $\sum_{k=0}^{t-1}\alpha_k$ and taking the limit superior as $t \to \infty$, it follows from Jensen's inequality that
$\mathbb{E}[f(\widehat{\mathbf{z}}_t)]\leq\frac{\sum_{k=0}^{t-1}\alpha_k\mathbb{E}[f(\bar{\mathbf{z}}_k)]}{\sum_{k=0}^{t-1}\alpha_k}$,
and Lemma~\ref{theorem:Stolz-Cesaro} that
$\frac{\sum_{k=0}^\infty\alpha_k\beta_{1,k}}{\sum_{k=0}^\infty\alpha_k} =\lim_{k\to\infty}\beta_{1,k}= 0$, $\frac{\sum_{k=0}^\infty\alpha^2_k}{\sum_{k=0}^\infty\alpha_k} = \alpha_\infty$,
we obtain
$\limsup_{k\to\infty}\mathbb{E}[f(\widehat{\mathbf{z}}_k)] - f^\star\leq \alpha_\infty\sum_{k=0}^\infty\tilde{\beta}_kN\mathcal{K}_2[\mathcal{K}_1+\Phi_1(\frac{(2N+\epsilon)\Gamma\gamma}{1-\gamma} + 2N)]+\alpha_\infty [\frac{(N\Phi_1(3\mathcal{K}_1+\bar{\beta}\mathcal{K}_2)+\Phi_2)\Gamma\gamma}{1-\gamma}+2.5\Phi_2+2\Phi_3+N\varsigma (\Phi_2+N(3\mathcal{K}_1+\bar{\beta}\mathcal{K}_2))\Gamma]$.

It follows from Assumption~\ref{assumption_local_f} and Proposition~\ref{prop:consensus}-(1) that
$\limsup_{k\to\infty}(\mathbb{E}[f(\widehat{\mathbf{x}}^i_k)] - \mathbb{E}[f(\widehat{\mathbf{z}}_k)])
\leq \hat{D}\limsup_{k\to\infty}\mathbb{E}[\|\widehat{\mathbf{x}}^i_k-\widehat{\mathbf{z}}_k\|]
\leq \hat{D}(N\varsigma\Gamma+\Phi_1+\frac{\Phi_1\Gamma\gamma}{1-\gamma})\alpha_\infty$.
The desired result follows by combining the preceding two relations.
\end{Proof}

\begin{Remark}
Theorem~\ref{theorem:non_summable} shows that the cost value of the multi-agent system will finally converge to a neighborhood of its optimal value with an error bounded by some terms, which are dependent on the step-size $\alpha_k$ and parameters $\beta_{1,k}, \beta_{2,k}$. Appropriate choice of the step-size and parameters will lead to the exact convergence to the optimal value. In particular, if the step-size $\alpha_k$ is set to $1/(k+1)^a$, where $a \in (0,1)$; the parameters $\beta_{1,k}, \beta_{2,k}$ are set to $1/(k+1)^{p_1}$ and $1/(k+1)^{p_2}$, respectively, where $p_1>0$ and $p_2-p_1>1$; then $\alpha_\infty = 0$ and $\sum_{k=0}^\infty\tilde{\beta}_k < \infty$, which means all the error terms will converge to 0. 
On the other hand, Theorem~\ref{theorem:non_summable} only proves the convergence of $\mathbb{E}[f(\widehat{\mathbf{x}}^i_k)]$, but cannot state anything about the convergence of the sequence $\widehat{\mathbf{x}}^i_k$, for $i\in\mathcal{V}$. We remark that achieving the exact convergence to the optimal value (\textit{i.e.}, $f(\widehat{\mathbf{x}}^i_k)\to f^\star$) is theoretically weaker than the exact convergence to an optimal solution (\textit{i.e.}, $\widehat{\mathbf{x}}^i_k \to x^\star$). 
The exact convergence of the sequence $\widehat{\mathbf{x}}^i_k$ to the optimal solution can be guaranteed based on the square-summable step-size condition, by using the Robbins-Siegmund's Lemma \cite{Polyak1987}, see the proof of Theorem~\ref{theorem:non_summable_square_summable}.
\end{Remark}

In the following corollary, we characterize the convergence rate of the proposed algorithm for both a diminishing step-size of $\alpha_k=\frac{\alpha}{\sqrt{k+2}}$ and a constant step-size of $\alpha_k=\frac{\alpha}{\sqrt{t+2}}$ if the number of iterations $t$ is known in advance.

\begin{Corollary}\label{corollary:convergence_rate}
Suppose Assumptions~\ref{assumption_graph}, \ref{assumption_local_f} and \ref{assumption_random_variables} hold. Let $\{\widehat{\mathbf{x}}^i_k\}_{k\geq0}$ be the sequence generated by \eqref{eq:algorithm} with a step-size sequence $\alpha_k$. Let $\epsilon$ be the constant such that $0<\epsilon <\min(\bar{\epsilon},\frac{1-\gamma}{2\sqrt{3}N\Gamma\gamma})$, where $\Gamma$ and $\gamma$ are some constants, and $\bar{\epsilon} = (\frac{1-|\lambda_3|}{20+8N})^N$ with $\lambda_3$ being the third largest eigenvalue of the weighting matrix $A$ in \eqref{eq:d-dgd} by setting $\epsilon =0$. Let the parameters $\beta_{1,k}, \beta_{2,k}$ be set to $\frac1{(k+2)^{p_1}}$ and $\frac1{(k+2)^{p_2}}$, respectively, where $p_1>1$ and $p = p_2-p_1>1$. Then 
\begin{enumerate}
\item if the step-size $\alpha_k=\frac{\alpha}{\sqrt{k+2}}$, $k=0,\ldots,t-1$, we have
\begin{align*}
\mathbb{E}[f(\widehat{\mathbf{x}}^i_t)] - f^\star \leq O(\ln t/\sqrt{t}),
\end{align*}
\item if the step-size $\alpha_k=\frac{\alpha}{\sqrt{t+2}}$, $k=0,\ldots,t-1$, we have
\begin{align*}
\mathbb{E}[f(\widehat{\mathbf{x}}^i_t)] - f^\star \leq O(1/\sqrt{t}).
\end{align*}
\end{enumerate}
\end{Corollary}

\begin{Proof}
Following the proof of Theorem~\ref{theorem:non_summable}, we can obtain that
\begin{align*}
&\mathbb{E}[f(\widehat{\mathbf{x}}^i_t)]- f^\star\\
&\quad\leq\frac1{\sum_{k=0}^{t-1}\alpha_k}\bigg[{C_0}+{C_1\sum_{k=0}^{t-1}\alpha^2_k}+{C_2\sum_{k=0}^{t-1}\alpha_k\beta_{1,k}}\\
&\quad\quad+{C_3\sum_{k=0}^{t-1}\alpha_k\tilde{\beta}_k}+{C_4\bigg(\sum_{k=0}^{t-1}\alpha^2_k\bigg)\bigg(\sum_{k=0}^{t-1}\tilde{\beta}_k\bigg)}\bigg],
\end{align*}
where $C_0, C_1, C_2, C_3$ and $C_4$ some constants.

For (1), $\alpha_k=\frac{\alpha}{\sqrt{k+2}}$, $k=0,\ldots,t-1$, we have
\begin{align*}
&\mathbb{E}[f(\widehat{\mathbf{x}}^i_t)] - f^\star \leq \frac{C_0}{2\alpha[\sqrt{t+2}-\sqrt{2}]}+ \frac{\alpha C_1\ln{(t+1)}}{2(\sqrt{t+2}-\sqrt{2})}\\
&\quad+ \frac{C_2(1-\frac1{(t+1)^{p_1-0.5}})}{[\sqrt{t+2}-\sqrt{2}](2p_1-1)}+ \frac{C_3(1-\frac1{(t+1)^{p-0.5}})}{[\sqrt{t+2}-\sqrt{2}](2p-1)}\\
&\quad+\frac{\alpha C_4\ln{(t+1)}(1-\frac1{(t+1)^{p-1}})}{2(\sqrt{t+2}-\sqrt{2})(p-1)}\\
&= O(1/\sqrt{t}) + O(\ln{t}/\sqrt{t})=O(\ln{t}/\sqrt{t}).
\end{align*}
Likewise for (2), $\alpha_k=\frac{\alpha}{\sqrt{t+2}}$, $k=0,\ldots,t-1$, we have
\begin{align*}
&\mathbb{E}[f(\widehat{\mathbf{x}}^i_t)] - f^\star \leq \frac{C_0\sqrt{t+2}}{t\alpha}+ \frac{\alpha C_1}{\sqrt{t+2}}\\
&\quad+ \frac{C_2(1-\frac1{(t+1)^{p_1-1}})}{t(p_1-1)}+ \frac{C_3(1-\frac1{(t+1)^{p-1}})}{t(p-1)}\\
&\quad+\frac{\alpha C_4(1-\frac1{(t+1)^{p-1}})}{\sqrt{t+2}(p-1)}\\
&= O(1/\sqrt{t}) + O(1/t)=O(1/\sqrt{t}).
\end{align*}
which gives the desired convergence rate results.
\end{Proof}

\section{Numerical Simulation}\label{sec:simulation}

In this section, we investigate the performance of the proposed algorithm through a numerical example. In particular, we consider a non-smooth test problem in a multi-agent system with $N$ agents originated from \cite{Nesterov2017}:
\begin{equation*}
\min f(\mathbf{x}) = \sum_{i=1}^{N} \bigg(l_i|x_1 - 1|+\sum_{d=1}^{n-1}|1+ x_{d+1}-2x_d|^2\bigg),\mathbf{x} \in \mathcal{X},
\end{equation*}
where $\mathbf{x} = [x_1, \ldots, x_n]^\top\in\mathcal{X}\subseteq\mathbb{R}^n$, $l_i, i = 1,2,\ldots,N$ is a positive constant.

In the simulation, the performance of the proposed algorithm is investigated from the following perspectives: the step-size and parameters selections, and comparison with both state-of-the-art gradient-free algorithm and gradient-based algorithm.
Throughout the simulation, we let $[A_r]_{ij} = 1/|\mathcal{N}^{\text{in}}_i|$ and $[A_c]_{ij} = 1/|\mathcal{N}^{\text{out}}_j|$, where $|\mathcal{N}|$ denotes the number of elements in $\mathcal{N}$. $l_i$ is randomly set in $[0.5,1.5]$.

\subsection{Influence of Step-Size $\alpha_k$ and Parameters $\beta_{1,k}, \beta_{2,k}$}

\begin{figure}[!t]
    \centering
    \includegraphics[width=2.3in]{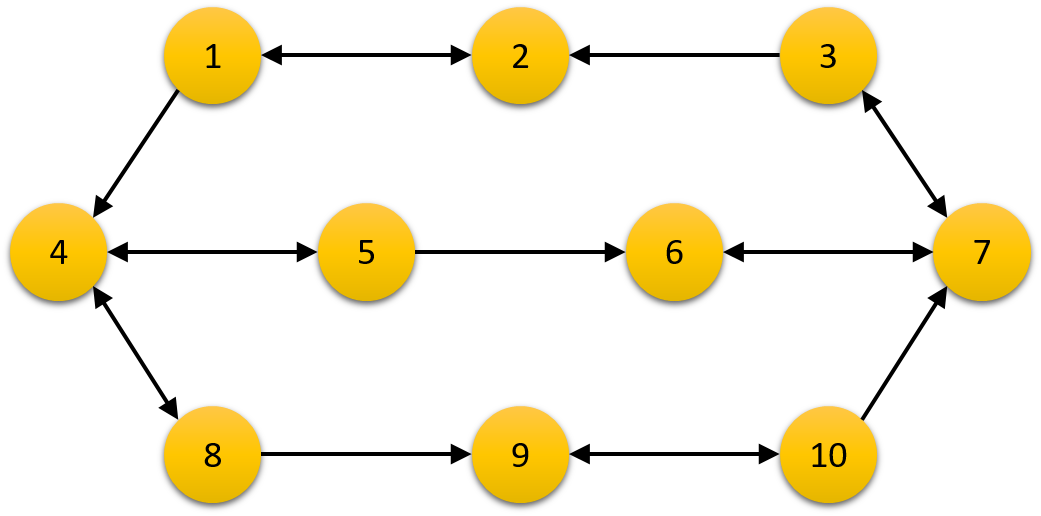}%
    \caption{Communication topology.}%
    \label{fig:network.PNG}%
\end{figure}

In this part, we set the dimension of the problem $n = 1$, the number of agents $N = 10$ under the directed graph $\mathcal{G}$ shown in Fig.~\ref{fig:network.PNG}. Then, we investigated the performance of the algorithm for the cases of different step-size $\alpha_k$ and two positive parameter sequences $\beta_{1,k}, \beta_{2,k}$, respectively. 

\begin{figure}[!t]
    \centering
    \includegraphics[width=3.2in]{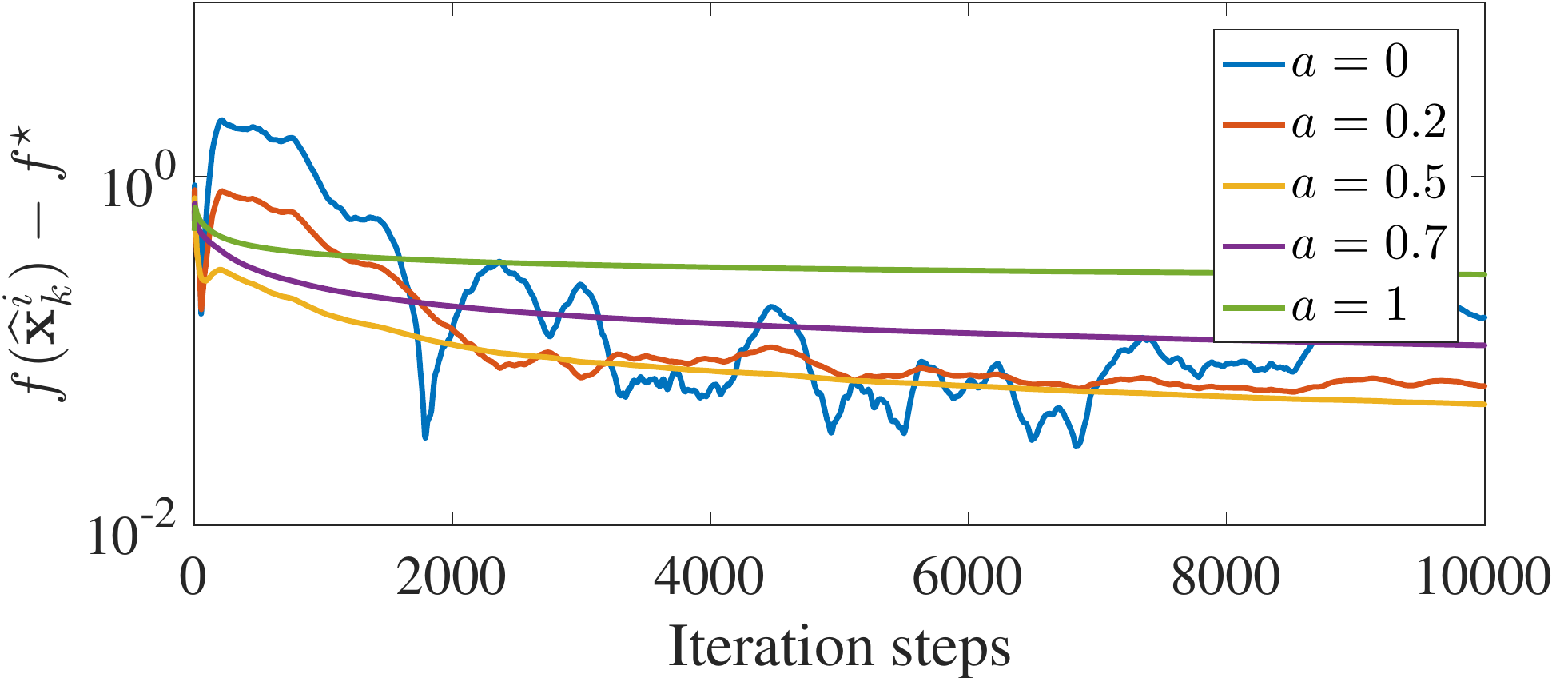}%
    \caption{Influence of step-size $\alpha_k$ on the convergence property.}%
    \label{fig:stepSize_compare.PNG}%
\end{figure}

To test the influence of the step-size on the convergence, we set the step-size $\alpha_k = 0.1/(1+k)^a$, where $a = 0, 0.2, 0.5, 0.7$ and $1$. It should be noted that the step-size $\alpha_k$ is not square-summable for $a = 0, 0.2, 0.5$. Two positive sequences were set to $\beta_{1,k} = 1/(1+k)^{1.5}$ and $\beta_{2,k} = 1/(1+k)^{2.5}$. The convergence result was shown in Fig.~\ref{fig:stepSize_compare.PNG}. As can be seen, both the optimality gap decreases for diminishing step-sizes, which is consistent with our findings in Theorem~\ref{theorem:non_summable}. Moreover, it can be observed that faster convergence result is attained with slower diminishing step-size (\textit{i.e.}, smaller $a$), but larger errors (oscillations in the plot) are incurred. 

\begin{figure}[!t]
    \centering
    \includegraphics[width=3.2in]{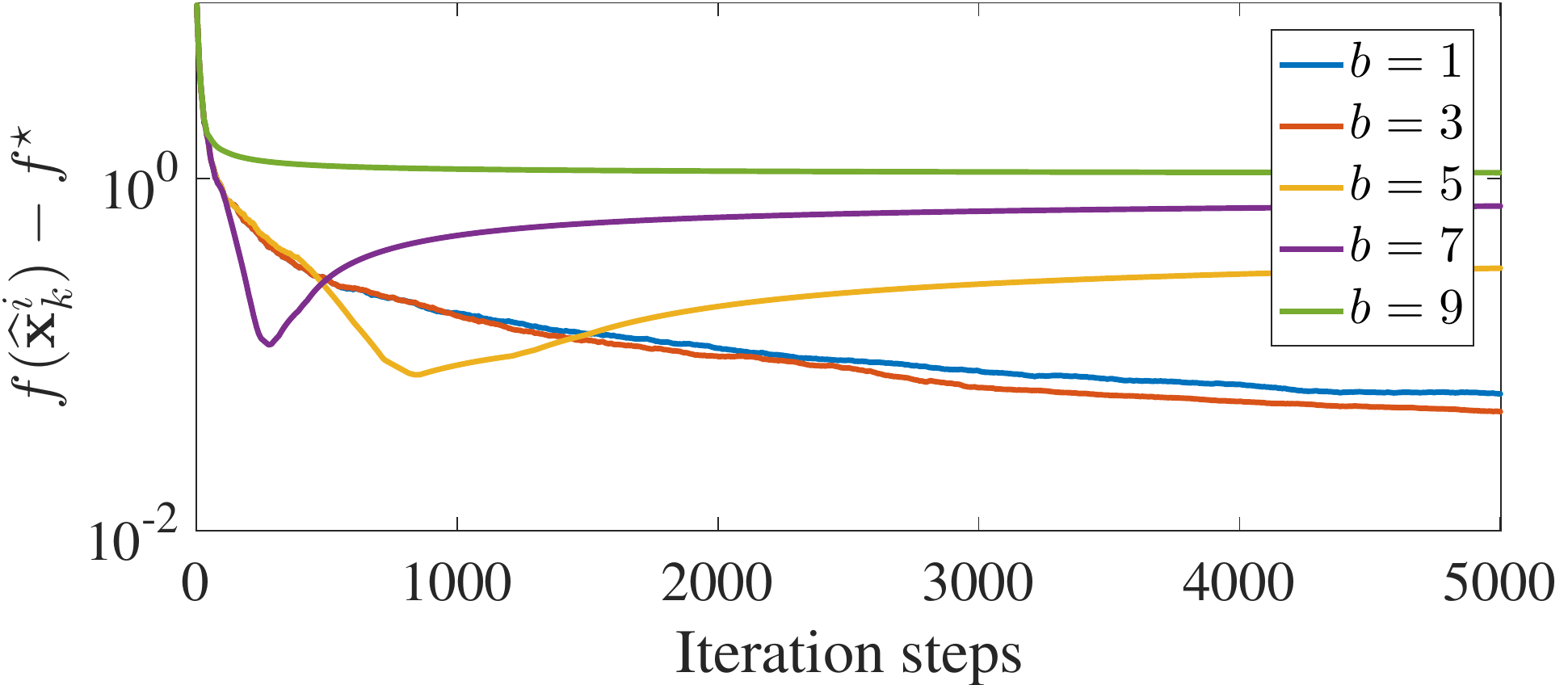}%
    \caption{Influence of $\tilde{\beta}_k$ on the convergence property.}%
    \label{fig:beta_compare.PNG}%
\end{figure}

To test the influence of the two positive parameter sequences on the convergence, we set $\beta_{1,k} = 1/(1+k)^{1.5}$, $\tilde{\beta}_k = \beta_{2,k}/\beta_{1,k} = 1/(1+k)^b$, where $b = 1,3,5,7$ and $9$. The step-size $\alpha_k$ was set to $0.1/\sqrt{k+1}$. The convergence result under these five cases was plotted in Fig.~\ref{fig:beta_compare.PNG}. As can be seen, typical $b$ values (ranging from 1 to 3) do not have much influence on the convergence rate. However, it can also be observed that when $b$ is increasing, the convergence performance is downgraded.

\subsection{Comparison with the State-Of-The-Art Algorithms}
\begin{figure}[!t]
    \centering
    \includegraphics[width=3.0in]{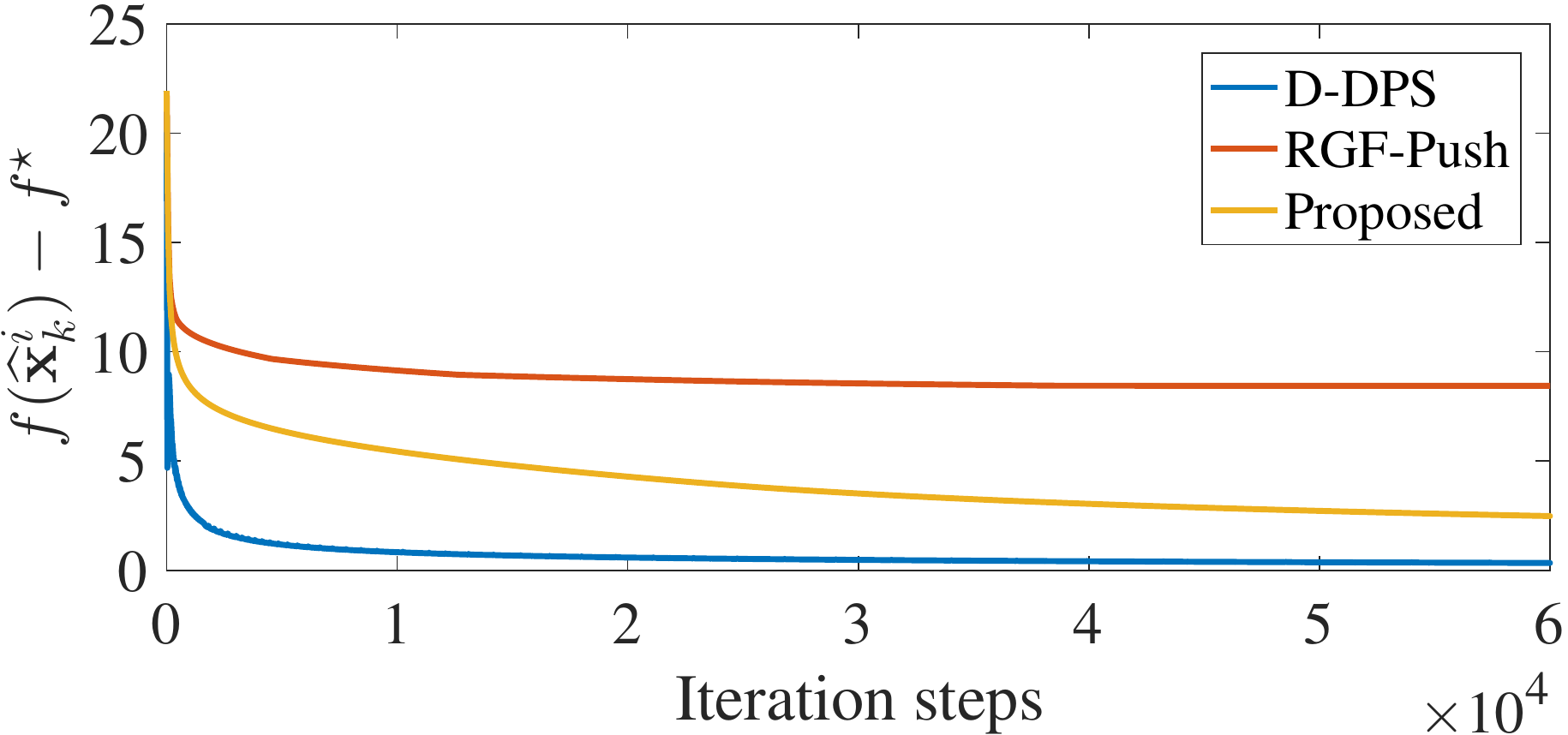}%
    \caption{Comparison between D-DPS, RGF-Push and the proposed method.}%
    \label{fig:algo_compare.PNG}%
\end{figure}

In this part, we compared our proposed method with the state-of-the-art algorithms, including the randomized gradient-free push-sum protocol (RGF-Push) proposed in \cite{Yuan2015a} using diminishing smoothing parameter and a subgradient-based method (D-DPS) proposed in \cite{Xi2016}. All these three methods can work for directed graphs. We set the dimension of the problem $n = 2$, the number of agents $N = 10$ under the directed graph $\mathcal{G}$ shown in Fig.~\ref{fig:network.PNG}. The step-size was set to $\alpha_k = 0.1/(k+1)^{0.5}$. The convergence results of all three methods were shown in Fig.~\ref{fig:algo_compare.PNG}. As can be seen, our proposed method shows a similar performance to the RGF-Push protocol, where both methods exhibit a theoretical convergence rate of $\ln k/\sqrt{k}$. The gradient-based algorithm (D-DPS) outperforms the two gradient-free methods as expected due to the use of the true gradient information. 

\section{Conclusions}\label{sec:conclusion}
This paper has considered a set constrained distributed optimization problem with possibly non-smooth cost functions. A distributed projected pseudo-gradient descent algorithm with an optimal averaging scheme has been proposed to solve the problem. 
The proposed algorithm has been shown to achieve the exact convergence to the optimal value with any positive, non-summable and non-increasing step-size sequence. When the step-size is also square-summable, the exact convergence to an optimal solution has been guaranteed.
Theoretical analysis on the convergence rate of the proposed algorithm has also been provided. To illustrate its performance, the proposed algorithm has been tested in a non-smooth problem. 
The convergence properties have been investigated, and the effectiveness has been verified by comparing with the state-of-the-art algorithms.

\section*{APPENDIX}

\subsection{Proof of Lemma~\ref{lemma:bound_on_g}}

For part (1), by definition of $g^i_k$ in \eqref{eq:d-dgd}
$\|g^i_k\|\leq \|\epsilon \mathbf{y}^i_k\|+\|\mathbf{x}^i_{k+1}-\sum_{j=1}^N [A_r]_{ij}\mathbf{x}^j_k\|
\leq \epsilon \|\mathbf{y}^i_k\|+ \|\epsilon \mathbf{y}^i_k - \alpha_k \mathbf{g}^i(\mathbf{x}^i_k)\|
\leq2\epsilon \|\mathbf{y}^i_k\| + \alpha_k\|\mathbf{g}^i(\mathbf{x}^i_k)\|$,
where the second inequality follows from the projection's nonexpansive property.
Summing over $i=1,\ldots,N$, and applying Lemma~\ref{lemma:bound_on_consensus}-(2),
\begin{align}
&\bm{G}_k \leq4N^2\varsigma\epsilon\Gamma\gamma^k\nonumber\\
&\quad\quad+\alpha_k\sum_{i=1}^N\|\mathbf{g}^i(\mathbf{x}^i_k)\|+ 2N\epsilon\Gamma\sum_{r=1}^{k-1}\gamma^{k-r}\bm{G}_{r-1}.\label{eq:sum_of_g}
\end{align}
Multiplying both sides by $\alpha_k$, summing over from $k=1$ to $K$, and noting that 
$\sum_{k=1}^K\alpha_k\sum_{r=1}^{k-1}\gamma^{k-r}\bm{G}_{r-1}\leq \sum_{k=1}^K\sum_{r=1}^{k-1}\gamma^{k-r}\alpha_r\bm{G}_{r-1}\leq\frac{\gamma}{1-\gamma}\sum_{k=1}^K\alpha_k\bm{G}_k$,
we obtain
\begin{align*}
&\sum_{k=1}^K\alpha_k\bm{G}_k \leq 2N^2\varsigma\epsilon\Gamma\sum_{k=1}^K\gamma^{2k}+2N^2\varsigma\epsilon\Gamma\sum_{k=1}^K\alpha_k^2\\
&\quad+\sum_{k=1}^K\alpha^2_k\sum_{i=1}^N\|\mathbf{g}^i(\mathbf{x}^i_k)\|+\frac{2N\epsilon\Gamma\gamma}{1-\gamma}\sum_{k=1}^K\alpha_k\bm{G}_k,
\end{align*}
Taking the total expectation and invoking Lemma~\ref{lemma:property_f_mu}-(3), we have
$\sum_{k=1}^K\alpha_k\mathbb{E}[\bm{G}_k] \leq 2N^2\varsigma\epsilon\Gamma\sum_{k=1}^K\gamma^{2k}
+N\mathcal{K}_1\sum_{k=1}^K\alpha_k^2+\frac{2N\epsilon\Gamma\gamma}{1-\gamma}\sum_{k=1}^K\alpha_k\mathbb{E}[\bm{G}_k]$.
Re-arranging the term and noticing that $\epsilon < \frac{1-\gamma}{2\sqrt{3}N\Gamma\gamma}< \frac{1-\gamma}{2N\Gamma\gamma}$, we obtain the desired result by denoting $\Phi_1 = \frac{N\mathcal{K}_1}{1-\gamma(2N\epsilon\Gamma+1)}$, and $\Psi_1 = \frac{2N^2\varsigma\epsilon\Gamma\gamma^2}{1-\gamma(2N\epsilon\Gamma+1)}$.

For part (2), squaring both sides of \eqref{eq:sum_of_g}, summing over from $k=1$ to $K$, and taking the total expectation, we have
\begin{align*}
&\sum_{k=1}^K\mathbb{E}[\bm{G}^2_k] \leq48N^4\varsigma^2\epsilon^2\Gamma^2\sum_{k=1}^K\gamma^{2k}+3N^2\mathcal{K}_1^2\sum_{k=1}^K\alpha_k^2\\
&\quad+ 12N^2\epsilon^2\Gamma^2\sum_{k=1}^K\mathbb{E}\bigg[\bigg(\sum_{r=1}^{k-1}\gamma^{k-r}\bm{G}_{r-1}\bigg)^2\bigg].
\end{align*}
Applying Cauchy-Schwarz inequality on the last term that
\begin{align*}
\bigg(\sum_{r=1}^{k-1}\gamma^{k-r}\bm{G}_{r-1}\bigg)^2\leq\frac{\gamma}{1-\gamma}\sum_{r=1}^{k-1}\gamma^{k-r}\bm{G}^2_{r-1},
\end{align*}
we obtain
\begin{align*}
\sum_{k=1}^K\mathbb{E}[\bm{G}^2_k] &\leq48N^4\varsigma^2\epsilon^2\Gamma^2\sum_{k=1}^K\gamma^{2k}+3N^2\mathcal{K}_1^2\sum_{k=1}^K\alpha_k^2\\
&\quad+\frac{12N^2\epsilon^2\Gamma^2\gamma}{1-\gamma}\sum_{k=1}^K\sum_{r=1}^{k-1}\gamma^{k-r}\mathbb{E}[\bm{G}^2_{r-1}]\\
&\leq48N^4\varsigma^2\epsilon^2\Gamma^2\sum_{k=1}^K\gamma^{2k}+3N^2\mathcal{K}_1^2\sum_{k=1}^K\alpha_k^2\\
&\quad+\frac{12N^2\epsilon^2\Gamma^2\gamma^2}{(1-\gamma)^2}\sum_{k=1}^K\mathbb{E}[\bm{G}_k^2].
\end{align*}
Re-arranging the term and noticing that $\epsilon < \frac{1-\gamma}{2\sqrt{3}N\Gamma\gamma}$, we obtain the desired result by denoting $\Phi_2 = \frac{3N^2\mathcal{K}_1^2}{(1-\gamma)^2-12N^2\epsilon^2\Gamma^2\gamma^2}$, and $\Psi_2 = \frac{48N^4\varsigma^2\epsilon^2\Gamma^2\gamma^2}{(1-\gamma)^2-12N^2\epsilon^2\Gamma^2\gamma^2}$.

For part (3), multiplying both sides of \eqref{eq:sum_of_g} by $\sum_{i=1}^N\alpha_k\|\mathbf{g}^i(\mathbf{x}^i_k)\|$, summing over from $k=1$ to $K$, and taking the total expectation, we have
\begin{align}
&\sum_{k=1}^K\sum_{i=1}^N\alpha_k\mathbb{E}[\|\mathbf{g}^i(\mathbf{x}^i_k)\|\bm{G}_k] \leq 2N^3\varsigma\epsilon\mathcal{K}_1\Gamma\sum_{k=1}^K\gamma^{2k}\nonumber\\
&\quad+2N^3\varsigma\epsilon\mathcal{K}_1\Gamma\sum_{k=1}^K\alpha_k^2+N^2\mathcal{K}_1^2\sum_{k=1}^K\alpha^2_k\nonumber\\
&\quad+2N\epsilon\Gamma\sum_{k=1}^K\sum_{r=1}^{k-1}\gamma^{k-r}\alpha_k\sum_{i=1}^N\mathbb{E}[\|\mathbf{g}^i(\mathbf{x}^i_k)\|\bm{G}_{r-1}].\label{eq:sum_alpha_g_G}
\end{align}
Based on Cauchy-Schwarz inequality that
\begin{align*}
\sum_{i=1}^N\mathbf{E}[\|\mathbf{g}^i(\mathbf{x}^i_k)\|\bm{G}_{r-1}] &\leq \sum_{i=1}^N\sqrt{\mathbf{E}[\|\mathbf{g}^i(\mathbf{x}^i_k)\|^2]\mathbf{E}[\bm{G}^2_{r-1}]} \\
&\leq N\mathcal{K}_1\sqrt{\mathbf{E}[\bm{G}^2_{r-1}]},
\end{align*}
the last term of \eqref{eq:sum_alpha_g_G} holds that
\begin{align*}
&2N\epsilon\Gamma\sum_{k=1}^K\sum_{r=1}^{k-1}\gamma^{k-r}\alpha_k\sum_{i=1}^N\mathbb{E}[\|\mathbf{g}^i(\mathbf{x}^i_k)\|\bm{G}_{r-1}]\\
&\leq2N^2\epsilon\mathcal{K}_1\Gamma\sum_{k=1}^K\sum_{r=1}^{k-1}\gamma^{k-r}\alpha_k\sqrt{\mathbf{E}[\bm{G}^2_{r-1}]}\\
&\leq N^2\epsilon\mathcal{K}_1\Gamma\sum_{k=1}^K\sum_{r=1}^{k-1}\gamma^{k-r}(\alpha_k^2+\mathbf{E}[\bm{G}^2_{r-1}])\\
&\leq \frac{N^2\epsilon\mathcal{K}_1\Gamma\gamma}{1-\gamma}\sum_{k=1}^K\alpha_k^2+\frac{N^2\epsilon\mathcal{K}_1\Gamma\gamma}{1-\gamma}\sum_{k=1}^K\mathbf{E}[\bm{G}^2_k]\\
&\leq \frac{N^2\epsilon\mathcal{K}_1\Gamma\gamma}{1-\gamma}\sum_{k=1}^K\alpha_k^2+\frac{N^2\epsilon\mathcal{K}_1\Gamma\gamma}{1-\gamma}\bigg(\Phi_2\sum_{k=1}^K\alpha_k^2+\Psi_2\bigg).
\end{align*}
Combining the above relation with \eqref{eq:sum_alpha_g_G}, we obtain the desired result by denoting $\Phi_3 = 2N^3\varsigma\epsilon\mathcal{K}_1\Gamma+N^2\mathcal{K}_1^2+\frac{N^2\epsilon\mathcal{K}_1\Gamma\gamma(1+\Phi_2)}{1-\gamma}$, and $\Psi_3 = \frac{2N^3\varsigma\epsilon\mathcal{K}_1\Gamma\gamma^2}{(1-\gamma)^2-12N^2\epsilon^2\Gamma^2\gamma^2}+\frac{N^2\epsilon\mathcal{K}_1\Gamma\gamma\Psi_2}{1-\gamma}$.



\bibliographystyle{IEEEtran}
\bibliography{d_dppgd_reference}

\end{document}